
\input amstex



\def\spaces{\space\space\space\space\space\space\space\space\space\space}
\def\spacess{\message{\spaces\spaces\spaces\spaces\spaces\spaces\spaces}}
\spacess
\spacess
\message{Annals of Mathematics Style: Current Version: 1.1. June 10, 1992}
\spacess
\spacess

\catcode`\@=11

\hyphenation{acad-e-my acad-e-mies af-ter-thought anom-aly anom-alies
an-ti-deriv-a-tive an-tin-o-my an-tin-o-mies apoth-e-o-ses apoth-e-o-sis
ap-pen-dix ar-che-typ-al as-sign-a-ble as-sist-ant-ship as-ymp-tot-ic
asyn-chro-nous at-trib-uted at-trib-ut-able bank-rupt bank-rupt-cy
bi-dif-fer-en-tial blue-print busier busiest cat-a-stroph-ic
cat-a-stroph-i-cally con-gress cross-hatched data-base de-fin-i-tive
de-riv-a-tive dis-trib-ute dri-ver dri-vers eco-nom-ics econ-o-mist
elit-ist equi-vari-ant ex-quis-ite ex-tra-or-di-nary flow-chart
for-mi-da-ble forth-right friv-o-lous ge-o-des-ic ge-o-det-ic geo-met-ric
griev-ance griev-ous griev-ous-ly hexa-dec-i-mal ho-lo-no-my ho-mo-thetic
ideals idio-syn-crasy in-fin-ite-ly in-fin-i-tes-i-mal ir-rev-o-ca-ble
key-stroke lam-en-ta-ble light-weight mal-a-prop-ism man-u-script
mar-gin-al meta-bol-ic me-tab-o-lism meta-lan-guage me-trop-o-lis
met-ro-pol-i-tan mi-nut-est mol-e-cule mono-chrome mono-pole mo-nop-oly
mono-spline mo-not-o-nous mul-ti-fac-eted mul-ti-plic-able non-euclid-ean
non-iso-mor-phic non-smooth par-a-digm par-a-bol-ic pa-rab-o-loid
pa-ram-e-trize para-mount pen-ta-gon phe-nom-e-non post-script pre-am-ble
pro-ce-dur-al pro-hib-i-tive pro-hib-i-tive-ly pseu-do-dif-fer-en-tial
pseu-do-fi-nite pseu-do-nym qua-drat-ics quad-ra-ture qua-si-smooth
qua-si-sta-tion-ary qua-si-tri-an-gu-lar quin-tes-sence quin-tes-sen-tial
re-arrange-ment rec-tan-gle ret-ri-bu-tion retro-fit retro-fit-ted
right-eous right-eous-ness ro-bot ro-bot-ics sched-ul-ing se-mes-ter
semi-def-i-nite semi-ho-mo-thet-ic set-up se-vere-ly side-step sov-er-eign
spe-cious spher-oid spher-oid-al star-tling star-tling-ly
sta-tis-tics sto-chas-tic straight-est strange-ness strat-a-gem strong-hold
sum-ma-ble symp-to-matic syn-chro-nous topo-graph-i-cal tra-vers-a-ble
tra-ver-sal tra-ver-sals treach-ery turn-around un-at-tached un-err-ing-ly
white-space wide-spread wing-spread wretch-ed wretch-ed-ly Brown-ian
Eng-lish Euler-ian Feb-ru-ary Gauss-ian Grothen-dieck Hamil-ton-ian
Her-mit-ian Jan-u-ary Japan-ese Kor-te-weg Le-gendre Lip-schitz
Lip-schitz-ian Mar-kov-ian Noe-ther-ian No-vem-ber Rie-mann-ian
Schwarz-schild Sep-tem-ber Za-mo-lod-chi-kov Knizh-nik quan-tum Op-dam
Mac-do-nald Ca-lo-ge-ro Su-ther-land Mo-ser Ol-sha-net-sky  Pe-re-lo-mov
in-de-pen-dent ope-ra-tors
}

\Invalid@\nofrills
\Invalid@\usualspace
\newif\ifnofrills@
\def\nofrills@#1#2{\relaxnext@
  \DN@{\ifx\next\nofrills
    \nofrills@true\let#2\relax\DN@\nofrills{\nextii@}%
  \else
    \nofrills@false\def#2{#1}\let\next@\nextii@\fi
\next@}}
\def\usualspace@#1{\ifnofrills@\def\usualspace{#1}\fi}
\def\addto#1#2{\csname \expandafter\eat@\string#1@\endcsname
  \expandafter{\the\csname \expandafter\eat@\string#1@\endcsname#2}}
\newdimen\bigsize@
\def\big@#1#2{{\hbox{$\left#2\vcenter to#1\bigsize@{}%
  \right.\nulldelimiterspace\z@\m@th$}}}
\def\big{\big@\@ne}
\def\Big{\big@{1.5}}
\def\bigg{\big@\tw@}
\def\Bigg{\big@{2.5}}
\def\raggedcenter@{\leftskip\z@ plus.4\hsize \rightskip\leftskip
 \parfillskip\z@ \parindent\z@ \spaceskip.3333em \xspaceskip.5em
 \pretolerance9999\tolerance9999 \exhyphenpenalty\@M
 \hyphenpenalty\@M \let\\\linebreak}
\def\upperspecialchars{\def\ss{SS}\let\i=I\let\j=J\let\ae\AE\let\oe\OE
  \let\o\O\let\aa\AA\let\l\L}
\def\uppercasetext@#1{%
  {\spaceskip1.2\fontdimen2\the\font plus1.2\fontdimen3\the\font
   \upperspecialchars\uctext@#1$\m@th\aftergroup\eat@$}}
\def\uctext@#1$#2${\endash@#1-\endash@$#2$\uctext@}
\def\endash@#1-#2\endash@{%
\uppercase{#1}\if\notempty{#2}--\endash@#2\endash@\fi}
\def\runaway@#1{\DN@{#1}\ifx\envir@\next@
  \Err@{You seem to have a missing or misspelled \string\end#1 ...}%
  \let\envir@\empty\fi}
\newif\iftemp@
\def\notempty#1{TT\fi\def\test@{#1}\ifx\test@\empty\temp@false
  \else\temp@true\fi \iftemp@}

\font@\tensmc=cmcsc10
\font@\sevenex=cmex7 
\font@\sevenit=cmti7
\font@\eightrm=cmr8 
\font@\sixrm=cmr6 
\font@\eighti=cmmi8     \skewchar\eighti='177 
\font@\sixi=cmmi6       \skewchar\sixi='177   
\font@\eightsy=cmsy8    \skewchar\eightsy='60 
\font@\sixsy=cmsy6      \skewchar\sixsy='60   
\font@\eightex=cmex8 %
\font@\eightbf=cmbx8 
\font@\sixbf=cmbx6   
\font@\eightit=cmti8 
\font@\eightsl=cmsl8 
\font@\eightsmc=cmcsc10 
\font@\eighttt=cmtt8 

\loadmsam
\loadmsbm
\loadeufm
\UseAMSsymbols

\def\penaltyandskip@#1#2{\relax\ifdim\lastskip<#2\relax\removelastskip
      \ifnum#1=\z@\else\penalty@#1\relax\fi\vskip#2%
  \else\ifnum#1=\z@\else\penalty@#1\relax\fi\fi}
\def\nobreak{\penalty\@M
  \ifvmode\def\penalty@{\let\penalty@\penalty\count@@@}%
  \everypar{\let\penalty@\penalty\everypar{}}\fi}
\let\penalty@\penalty

\def\block{\RIfMIfI@\nondmatherr@\block\fi
       \else\ifvmode\vskip\abovedisplayskip\noindent\fi
        $$\def\endblock{\par\egroup$$}\fi
  \vbox\bgroup\advance\hsize-2\indenti\noindent}
\def\endblock{\par\egroup}

\def\logo@{\baselineskip2pc \hbox to\hsize{\hfil\eightpoint Typeset by
 \AmSTeX}}




\font\elevensc=cmcsc10 scaled\magstephalf
\font\tensc=cmcsc10

\font\eightsc=cmcsc10 scaled800

\font\elevenrm=cmr10 scaled \magstephalf
\font\ninerm=cmr9
\font\eightrm=cmr8
\font\sixrm=cmr6
\font\fiverm=cmr5

\font\eleveni=cmmi10 scaled\magstephalf
\font\ninei=cmmi9
\font\eighti=cmmi8
\font\sixi=cmmi6
\font\fivei=cmmi5
\skewchar\ninei='177 \skewchar\eighti='177 \skewchar\sixi='177
\skewchar\eleveni='177

\font\elevensy=cmsy10 scaled\magstephalf
\font\ninesy=cmsy9
\font\eightsy=cmsy8
\font\sixsy=cmsy6
\font\fivesy=cmsy5
\skewchar\ninesy='60 \skewchar\eightsy='60 \skewchar\sixsy='60
\skewchar\elevensy'60

\font\eighteenbf=cmbx10 scaled\magstep3

\font\twelvebf=cmbx10 scaled \magstep1
\font\elevenbf=cmbx10 scaled \magstephalf
\font\tenbf=cmbx10
\font\ninebf=cmbx9
\font\eightbf=cmbx8
\font\sixbf=cmbx6
\font\fivebf=cmbx5

\font\elevenit=cmti10 scaled\magstephalf
\font\nineit=cmti9
\font\eightit=cmti8

\font\eighteenmib=cmmib10 scaled \magstep3
\font\twelvemib=cmmib10 scaled \magstep1
\font\elevenmib=cmmib10 scaled\magstephalf
\font\tenmib=cmmib10
\font\eightmib=cmmib10 scaled 800 
\font\sixmib=cmmib10 scaled 600

\font\eighteensyb=cmbsy10 scaled \magstep3 
\font\twelvesyb=cmbsy10 scaled \magstep1
\font\elevensyb=cmbsy10 scaled \magstephalf
\font\tensyb=cmbsy10 
\font\eightsyb=cmbsy10 scaled 800
\font\sixsyb=cmbsy10 scaled 600
 
\font\elevenex=cmex10 scaled \magstephalf
\font\tenex=cmex10     
\font\eighteenex=cmex10 scaled \magstep3


\def\elevenpoint{\def\rm{\fam0\elevenrm}%
  \textfont0=\elevenrm \scriptfont0=\eightrm \scriptscriptfont0=\sixrm
  \textfont1=\eleveni \scriptfont1=\eighti \scriptscriptfont1=\sixi
  \textfont2=\elevensy \scriptfont2=\eightsy \scriptscriptfont2=\sixsy
  \textfont3=\elevenex \scriptfont3=\tenex \scriptscriptfont3=\tenex
  \def\bf{\fam\bffam\elevenbf}%
  \def\it{\fam\itfam\elevenit}%
  \textfont\bffam=\elevenbf \scriptfont\bffam=\eightbf
   \scriptscriptfont\bffam=\sixbf
\normalbaselineskip=13.95pt
  \setbox\strutbox=\hbox{\vrule height9.5pt depth4.4pt width0pt\relax}%
  \normalbaselines\rm}

\elevenpoint 

\def\ninepoint{\def\rm{\fam0\ninerm}%
  \textfont0=\ninerm \scriptfont0=\sixrm \scriptscriptfont0=\fiverm
  \textfont1=\ninei \scriptfont1=\sixi \scriptscriptfont1=\fivei
  \textfont2=\ninesy \scriptfont2=\sixsy \scriptscriptfont2=\fivesy
  \textfont3=\tenex \scriptfont3=\tenex \scriptscriptfont3=\tenex
  \def\it{\fam\itfam\nineit}%
  \textfont\itfam=\nineit
  \def\bf{\fam\bffam\ninebf}%
  \textfont\bffam=\ninebf \scriptfont\bffam=\sixbf
   \scriptscriptfont\bffam=\fivebf
\normalbaselineskip=11pt
  \setbox\strutbox=\hbox{\vrule height8pt depth3pt width0pt\relax}%
  \normalbaselines\rm}

\def\eightpoint{\def\rm{\fam0\eightrm}%
  \textfont0=\eightrm \scriptfont0=\sixrm \scriptscriptfont0=\fiverm
  \textfont1=\eighti \scriptfont1=\sixi \scriptscriptfont1=\fivei
  \textfont2=\eightsy \scriptfont2=\sixsy \scriptscriptfont2=\fivesy
  \textfont3=\tenex \scriptfont3=\tenex \scriptscriptfont3=\tenex
  \def\it{\fam\itfam\eightit}%
  \textfont\itfam=\eightit
  \def\bf{\fam\bffam\eightbf}%
  \textfont\bffam=\eightbf \scriptfont\bffam=\sixbf
   \scriptscriptfont\bffam=\fivebf
\normalbaselineskip=12pt
  \setbox\strutbox=\hbox{\vrule height8.5pt depth3.5pt width0pt\relax}%
  \normalbaselines\rm}


\def\eighteenbold{\def\rm{\fam0\eighteenbf}%
  \textfont0=\eighteenbf \scriptfont0=\twelvebf \scriptscriptfont0=\tenbf
  \textfont1=\eighteenmib \scriptfont1=\twelvemib\scriptscriptfont1=\tenmib
  \textfont2=\eighteensyb \scriptfont2=\twelvesyb\scriptscriptfont2=\tensyb
  \textfont3=\eighteenex \scriptfont3=\tenex \scriptscriptfont3=\tenex
  \def\bf{\fam\bffam\eighteenbf}%
  \textfont\bffam=\eighteenbf \scriptfont\bffam=\twelvebf
   \scriptscriptfont\bffam=\tenbf
\normalbaselineskip=20pt
  \setbox\strutbox=\hbox{\vrule height13.5pt depth6.5pt width0pt\relax}%
\everymath {\fam0 }
\everydisplay {\fam0 }
  \normalbaselines\rm}

\def\elevenbold{\def\rm{\fam0\elevenbf}%
  \textfont0=\elevenbf \scriptfont0=\eightbf \scriptscriptfont0=\sixbf
  \textfont1=\elevenmib \scriptfont1=\eightmib \scriptscriptfont1=\sixmib
  \textfont2=\elevensyb \scriptfont2=\eightsyb \scriptscriptfont2=\sixsyb
  \textfont3=\elevenex \scriptfont3=\elevenex \scriptscriptfont3=\elevenex
  \def\bf{\fam\bffam\elevenbf}%
  \textfont\bffam=\elevenbf \scriptfont\bffam=\eightbf
   \scriptscriptfont\bffam=\sixbf
\normalbaselineskip=14pt
  \setbox\strutbox=\hbox{\vrule height10pt depth4pt width0pt\relax}%
\everymath {\fam0 }
\everydisplay {\fam0 }
  \normalbaselines\bf}

\hsize=31pc
\vsize=48pc

\parindent=22pt
\parskip=0pt

\widowpenalty=10000
\clubpenalty=10000

\topskip=12pt 

\skip\footins=20pt
\dimen\footins=3in 

\abovedisplayskip=6.95pt plus3.5pt minus 3pt
\belowdisplayskip=\abovedisplayskip


\voffset=7pt\hoffset= .7in

\newif\iftitle

\def\amheadline{\iftitle%
\hbox to\hsize{\hss\currannalsline\hss}\else\line{\ifodd\pageno
\hfill\thetitle\hfill\llap{\elevenrm\folio}\else\rlap{\elevenrm\folio}
\hfill\theauthors\hfill\fi}\fi}

\headline={\amheadline}
\footline={\global\titlefalse}


\def\annalsline#1#2{\vfill\eject
\ifodd\pageno\else 
\line{\hfill}
\vfill\eject\fi
\global\titletrue
\def\currannalsline{\eightrm 
{\eightbf#1} (#2), \thepages}}

\def\titleheadline#1{\def\one{#1}\ifx\one\empty\else
\def\thetitle{{
\let\\ \relax\eightsc\uppercase{#1}}}\fi}

\newif\ifshort

\let\shorttitle\titleheadline

\def\onpages#1#2{\def\thepages{#1--#2}}

\def\thismuchskip[#1]{\vskip#1pt}
\def\ilook{\ifx\next[ \let\go\thismuchskip\else
\let\go\relax\vskip1pt\fi\go}

\def\institution#1{\def\theinstitutions{\vbox{\baselineskip10pt
\def\and{{\eightrm and }}
\def\\{\futurelet\next\ilook}\eightsc #1}}}
\let\institutions\institution

\newwrite\auxfile

\def\startingpage#1{\def\one{#1}\ifx\one\empty\global\pageno=1\else
\global\pageno=#1\fi
\theoremcount=0 \eqcount=0 \sectioncount=0 
\openin1 \jobname.aux \ifeof1 
\onpages{#1}{???}
\else\closein1 \relax\input \jobname.aux
\onpages{#1}{\lastpage}
\fi\immediate\openout\auxfile=\jobname.aux
}

\def\endarticle{\ifRefsUsed\global\RefsUsedfalse%
\else\vskip21pt\theinstitutions%
\nobreak\vskip8pt
\fi%
\write\auxfile{\string\def\string\lastpage{\the\pageno}}}

\outer\def\bye{\endarticle\par \vfill \supereject \end}

\def\document{\let\fontlist@\relax\let\alloclist@\relax
 \elevenpoint}


\newif\ifacks
\long\def\acknowledgements#1{\def\one{#1}\ifx\one\empty\else
\vskip-\baselineskip%
\global\ackstrue\footnote{\ \unskip}{*#1}\fi}

\def\title#1{\titleheadline{#1}
\vbox to80pt{\vfill
\baselineskip=18pt
\parindent=0pt
\overfullrule=0pt
\hyphenpenalty=10000
\everypar={\hskip\parfillskip\relax}
\hbadness=10000
\def\\ {\vskip1sp}
\eighteenbold#1\vskip1sp}}

\newif\ifauthor

\def\author#1{\vskip11pt
\hbox to\hsize{\hss\tenrm By \tensc#1\ifacks\global\acksfalse*\fi\hss}
\ifshort\else\xdef\theauthors{{\eightsc\uppercase{#1}}}\fi%
\vskip21pt\global\authortrue\everypar={\global\authorfalse\everypar={}}}

\def\twoauthors#1#2{\vskip11pt
\hbox to\hsize{\hss%
\tenrm By \tensc#1 {\tenrm and} #2\ifacks\global\acksfalse*\fi\hss}
\ifshort\else\xdef\theauthors{{\eightsc\uppercase{#1 and #2}}}\fi%
\vskip21pt
\global\authortrue\everypar={\global\authorfalse\everypar={}}}


\newcount\theoremcount
\newcount\sectioncount
\newcount\eqcount

\newif\ifspecialnumon

\def\eqnumber=#1 {\global\eqcount=#1 \global\advance\eqcount by-1\relax}
\def\sectionnumber=#1 {\global\sectioncount=#1 
\global\advance\sectioncount by-1\relax}
\def\proclaimnumber=#1 {\global\theoremcount=#1 
\global\advance\theoremcount by-1\relax}

\newif\ifsection
\newif\ifsubsection

\def\intro{\global\authorfalse
\centerline{\bf Introduction}\everypar={}\vskip6pt}

\def\elevenboldmath#1{$#1$\egroup}
\def\mathbold{\hbox\bgroup\elevenbold\elevenboldmath}

\def\section#1{\global\theoremcount=0
\global\eqcount=0
\ifauthor\global\authorfalse\else%
\vskip18pt plus 18pt minus 6pt\fi%
{\parindent=0pt
\everypar={\hskip\parfillskip}
\def\\ {\vskip1sp}\elevenpoint\bf%
\ifspecialnumon\global\specialnumonfalse$\rm\spnum$%
\gdef\sectnum{$\rm\spnum$}%
\else\interlinepenalty=10000%
\global\advance\sectioncount by1\relax\the\sectioncount%
\gdef\sectnum{\the\sectioncount}%
\fi. \hskip6pt#1
\vrule width0pt depth12pt}
\hskip\parfillskip
\global\sectiontrue%
\everypar={\global\sectionfalse\global\interlinepenalty=0\everypar={}}%
\ignorespaces

}


\newif\ifspequation

\let\eqno\leqno 

\newif\ifineqalignno
\let\saveleqalignno\leqalignno                        
\def\leqalignno{\let\eqnu\Eeqnu\saveleqalignno}

\let\eqalignno\leqalignno

\def\sectandeqnum{%
\ifspecialnumon\global\specialnumonfalse
$\rm\spnum$\gdef\eqnum{{$\rm\spnum$}}\else\global\firstlettertrue
\global\advance\eqcount by1 
\ifappend\applett\else\the\sectioncount\fi.%
\the\eqcount
\xdef\eqnum{\ifappend\applett\else\the\sectioncount\fi.\the\eqcount}\fi}

\def\eqnu{\leqno{\hbox{\elevenrm\ifspequation\else(\fi\sectandeqnum
\ifspequation\global\spequationfalse\else)\fi}}}      

\def\Speqnu{\global\setbox\leqnobox=\hbox{\elevenrm
\ifspequation\else%
(\fi\sectandeqnum\ifspequation\global\spequationfalse\else)\fi}}

\def\Eeqnu{\hbox{\elevenrm
\ifspequation\else%
(\fi\sectandeqnum\ifspequation\global\spequationfalse\else)\fi}}

\newif\iffirstletter
\global\firstlettertrue
\def\eqletter#1{\global\specialnumontrue\iffirstletter\global\firstletterfalse
\global\advance\eqcount by1\fi
\gdef\spnum{\the\sectioncount.\the\eqcount#1}\eqnu}

\newbox\leqnobox
\def\outsideeqnu#1{\global\setbox\leqnobox=\hbox{#1}}

\def\eatone#1{}

\def\dosplit#1#2{\vskip-.5\abovedisplayskip
\setbox0=\hbox{$\let\eqno\outsideeqnu%
\let\eqnu\Speqnu\let\leqno\outsideeqnu#2$}%
\setbox1\vbox{\noindent\hskip\wd\leqnobox\ifdim\wd\leqnobox>0pt\hskip1em\fi%
$\displaystyle#1\mathstrut$\hskip0pt plus1fill\relax
\vskip1pt
\line{\hfill$\let\eqnu\eatone\let\leqno\eatone%
\displaystyle#2\mathstrut$\ifmathqed~~\qed\fi}}%
\copy1
\ifvoid\leqnobox
\else\dimen0=\ht1 \advance\dimen0 by\dp1
\vskip-\dimen0
\vbox to\dimen0{\vfill
\hbox{\unhbox\leqnobox}
\vfill}
\fi}

\everydisplay{\lookforbreak}

\long\def\lookforbreak #1$${\def\mathone{#1}
\expandafter\testforbreak\mathone\splitmath @}

\def\testforbreak#1\splitmath #2@{\def\mathtwo{#2}\ifx\mathtwo\empty%
#1$$%
\ifmathqed\vskip-\belowdisplayskip
\setbox0=\vbox{\let\eqno\relax\let\eqnu\relax$\displaystyle#1$}%
\vskip-\ht0\vskip-3.5pt\hbox to\hsize{\hfill\qed}
\vskip\ht0\vskip3.5pt\fi
\else$$\vskip-\belowdisplayskip
\vbox{\dosplit{#1}{\let\eqno\eatone
\let\splitmath\relax#2}}%
\nobreak\vskip.5\belowdisplayskip
\noindent\ignorespaces\fi}


\newif\ifmathqed



\newcount\linenum
\newcount\colnum

\def\spline{\omit&\multispan{\the\colnum}{\hrulefill}\cr}
\def\colcounter{\ifnum\linenum=1\global\advance\colnum by1\fi}

\def\everyline{\noalign{\global\advance\linenum by1\relax}%
\ifnum\linenum=2\spline\fi}

\def\mtable{\bgroup\offinterlineskip
\everycr={\everyline}\global\linenum=0
\halign\bgroup\vrule height 10pt depth 4pt width0pt
\hfill$##$\hfill\hskip6pt\ifnum\linenum>1
\vrule\fi&&\colcounter\hskip12pt\hfill$##$\hfill\hskip12pt\cr}

\def\endmtable{\crcr\egroup\egroup}




\def\xast{*}
\newcount\intable
\newcount\mathcol
\newcount\savemathcol
\newcount\topmathcol
\newdimen\arrayhspace
\newdimen\arrayvspace

\arrayhspace=8pt 
\arrayvspace=12pt 

\newif\ifdollaron

\def\mathalign#1{\def\arg{#1}\ifx\arg\xast%
\let\go\relax\else\let\go\mathalign%
\global\advance\mathcol by1 %
\global\advance\topmathcol by1 %
\expandafter\def\csname  mathcol\the\mathcol\endcsname{#1}%
\fi\go}

\def\arraypickapart#1]#2*{\if#1c \ifmmode\vcenter\else
\global\dollarontrue$\vcenter\fi\else%
\if#1t\vtop\else\if#1b\vbox\fi\fi\fi\bgroup%
\def\one{#2}}

\def\arraystrut{\vrule height .7\arrayvspace depth .3\arrayvspace width 0pt}

\def\array#1#2*{\def\firstarg{#1}%
\if\firstarg[ \def\two{#2} \expandafter\arraypickapart\two*\else%
\ifmmode\vcenter\else\vbox\fi\bgroup \def\one{#1#2}\fi%
\global\everycr={\noalign{\global\mathcol=\savemathcol\relax}}%
\def\\ {\cr}%
\global\advance\intable by1 %
\ifnum\intable=1 \global\mathcol=0 \savemathcol=0 %
\else \global\advance\mathcol by1 \savemathcol=\mathcol\fi%
\expandafter\mathalign\one*%
\mathcol=\savemathcol %
\halign\bgroup&\hskip.5\arrayhspace\arraystrut%
\global\advance\mathcol by1 \relax%
\expandafter\if\csname mathcol\the\mathcol\endcsname r\hfill\else%
\expandafter\if\csname mathcol\the\mathcol\endcsname c\hfill\fi\fi%
$\displaystyle##$%
\expandafter\if\csname mathcol\the\mathcol\endcsname r\else\hfill\fi\relax%
\hskip.5\arrayhspace\cr}

\def\endarray{\crcr\egroup\egroup%
\global\mathcol=\savemathcol %
\global\advance\intable by -1\relax%
\ifnum\intable=0 %
\ifdollaron\global\dollaronfalse $\fi
\loop\ifnum\topmathcol>0 %
\expandafter\def\csname  mathcol\the\topmathcol\endcsname{}%
\global\advance\topmathcol by-1 \repeat%
\global\everycr={}\fi%
}

\def\big#1{{\hbox{$\left#1\vbox to 10pt{}\right.\n@space$}}}
\def\Big#1{{\hbox{$\left#1\vbox to 13pt{}\right.\n@space$}}}
\def\bigg#1{{\hbox{$\left#1\vbox to 16pt{}\right.\n@space$}}}
\def\Bigg#1{{\hbox{$\left#1\vbox to 19pt{}\right.\n@space$}}}


\def\figcaption#1#2#3{\topinsert
\vskip4pt 
\vbox to#3{\vfill}\vskip1sp
\setbox0=\hbox{\eightsc Figure #1.\hskip12pt\eightpoint #2}
\ifdim\wd0>\hsize
\noindent\eightsc Figure #1.\hskip12pt\eightpoint #2
\else
\centerline{\eightsc Figure #1.\hskip12pt\eightpoint #2}
\fi
\vskip16pt
\endinsert}

\def\wfig#1#2#3{\topinsert
\vskip4pt 
\hbox to\hsize{\hss\vbox{\hrule height .25pt width #3
\hbox to #3{\vrule width .25pt height #2\hfill\vrule width .25pt height #2}
\hrule height.25pt}\hss}
\vskip1sp
\centerline{\eightsc Figure #1}
\vskip16pt
\endinsert}

\def\wfigcaption#1#2#3#4{\topinsert
\vskip4pt 
\hbox to\hsize{\hss\vbox{\hrule height .25pt width #4
\hbox to #4{\vrule width .25pt height #3\hfill\vrule width .25pt height #3}
\hrule height.25pt}\hss}
\vskip1sp
\setbox0=\hbox{\eightsc Figure #1.\hskip12pt\eightpoint\rm #2}
\ifdim\wd0>\hsize
\noindent\eightsc Figure #1.\hskip12pt\eightpoint\rm #2\else
\centerline{\eightsc Figure #1.\hskip12pt\eightpoint\rm #2}\fi
\vskip16pt
\endinsert}

\def\tabcaption#1#2{\vskip6pt
\setbox0=\hbox{\eightsc Table #1.\hskip12pt\eightpoint #2}
\ifdim\wd0>\hsize
\noindent\eightsc Table #1.\hskip12pt\eightpoint #2
\else
\centerline{\eightsc Table #1.\hskip12pt\eightpoint #2}
\fi
\vskip6pt}

\def\endinsert{\egroup\if@mid\dimen@\ht\z@\advance\dimen@\dp\z@ 
\advance\dimen@ 12\p@\advance\dimen@\pagetotal\ifdim\dimen@ >\pagegoal 
\@midfalse\p@gefalse\fi\fi\if@mid\smallskip\box\z@\bigbreak\else
\insert\topins{\penalty 100 \splittopskip\z@skip\splitmaxdepth\maxdimen
\floatingpenalty\z@\ifp@ge\dimen@\dp\z@\vbox to\vsize {\unvbox \z@ 
\kern -\dimen@ }\else\box\z@\nobreak\smallskip\fi}\fi\endgroup}

\def\pagecontents{
\ifvoid\topins \else\iftitle\else 
\unvbox \topins \fi\fi \dimen@ =\dp \@cclv \unvbox 
\@cclv 
\ifvoid\topins\else\iftitle\unvbox\topins\fi\fi
\ifvoid \footins \else \vskip \skip \footins \footnoterule 
\unvbox \footins \fi \ifr@ggedbottom \kern -\dimen@ \vfil \fi}


\newif\ifappend

\def\appendix#1#2{\def\applett{#1}\def\two{#2}%
\global\appendtrue
\global\theoremcount=0
\global\eqcount=0
\vskip18pt plus 18pt
\vbox{\parindent=0pt
\everypar={\hskip\parfillskip}
\def\\ {\vskip1sp}\elevenbold Appendix%
\ifx\applett\empty\gdef\applett{A}\ifx\two\empty\else.\fi%
\else\ #1.\fi\hskip6pt#2\vskip12pt}%
\global\sectiontrue%
\everypar={\global\sectionfalse\everypar={}}\nobreak\ignorespaces}

\newif\ifRefsUsed
\long\def\references{\global\RefsUsedtrue\vskip21pt
\theinstitutions
\global\everypar={}\global\bibnum=0
\vskip20pt\goodbreak\bgroup
\vbox{\centerline{\eightsc References}\vskip6pt}%
\ifdim\maxbibwidth>0pt
\leftskip=\maxbibwidth%
\parindent=-\maxbibwidth%
\else
\leftskip=18pt%
\parindent=-18pt%
\fi
\ninepoint
\frenchspacing
\nobreak\ignorespaces\everypar={\amref}%
}

\def\endreferences{\vskip1sp\egroup\global\everypar={}%
\nobreak\vskip8pt\vbox{\thereceived\therevised}
}

\newcount\bibnum

\def\amref#1 {\global\advance\bibnum by1%
\immediate\write\auxfile{\string\expandafter\string\def\string\csname
\space #1croref\string\endcsname{[\the\bibnum]}}%
\leavevmode\hbox to18pt{\hbox to13.2pt{\hss[\the\bibnum]}\hfill}}

\def\bibline{\hbox to30pt{\hrulefill}\/\/}

\def\name#1{{\eightsc#1}}

\newdimen\maxbibwidth
\def\AuthorRefNames [#1] {%
\immediate\write\auxfile{\string\def\string\cite\string##1{[\string##1]}}

\def\amref{\spamref}
\setbox0=\hbox{[#1] }\global\maxbibwidth=\wd0\relax}

\def\spamref[#1] {\leavevmode\hbox to\maxbibwidth{\hss[#1]\hfill}}


\def\footnoterule{\kern-3pt\hrule width1in height.5pt\kern2.5pt}

\def\footnote#1#2{%
\plainfootnote{#1}{{\eightpoint\normalbaselineskip11pt
\normalbaselines#2}}}

\def\vfootnote#1{%
\insert \footins \bgroup \eightpoint\baselineskip11pt
\interlinepenalty \interfootnotelinepenalty
\splittopskip \ht \strutbox \splitmaxdepth \dp \strutbox \floatingpenalty 
\@MM \leftskip \z@skip \rightskip \z@skip \spaceskip \z@skip 
\xspaceskip \z@skip
{#1}$\,$\footstrut \futurelet \next \fo@t}


\newif\iffirstadded
\newif\ifadded

\def\addedlett{}

\def\alltheoremnums{%
\ifspecialnumon\global\specialnumonfalse
\ifadded\global\addedfalse
\iffirstadded\global\firstaddedfalse
\global\advance\theoremcount by1 \fi
\ifappend\applett\else\the\sectioncount\fi.\the\theoremcount\addedlett%
\xdef\theoremnum{\ifappend\applett\else\the\sectioncount\fi.%
\the\theoremcount\addedlett}%
\else$\rm\spnum$\def\theoremnum{{$\rm\spnum$}}\fi%
\else\global\firstaddedtrue
\global\advance\theoremcount by1 
\ifappend\applett\else\the\sectioncount\fi.\the\theoremcount%
\xdef\theoremnum{\ifappend\applett\else\the\sectioncount\fi.%
\the\theoremcount}\fi}

\def\allcorolnums{%
\ifspecialnumon\global\specialnumonfalse
\ifadded\global\addedfalse
\iffirstadded\global\firstaddedfalse
\global\advance\corolcount by1 \fi
\the\corolcount\addedlett%
\else$\rm\spnum$\def\corolnum{$\rm\spnum$}\fi%
\else\global\advance\corolcount by1 
\the\corolcount\fi}


\newcount\corolcount
\def\xcorol{Corollary}
\def\xtheorem{Theorem}
\def\xmaintheorem{Main Theorem}

\newif\ifthtitle

\let\saverparen)
\let\savelparen(
\def\rmparenl{{\rm(}}
\def\rmparenr{{\rm\/)}}
{
\catcode`(=13
\catcode`)=13
\gdef\makeparensRM{\catcode`(=13\catcode`)=13\let(=\rmparenl%
\let)=\rmparenr%
\everymath{\let(\savelparen%
\let)\saverparen}%
\everydisplay{\let(\savelparen%
\let)\saverparen\lookforbreak}}}

\medskipamount=8pt plus.1\baselineskip minus.05\baselineskip

\def\rmtext#1{\hbox{\rm#1}}

\def\proclaim#1{\vskip-\lastskip
\def\one{#1}\ifx\one\xtheorem\global\corolcount=0\fi
\ifsection\global\sectionfalse\vskip-6pt\fi
\medskip
{\elevensc#1}%
\ifx\one\xmaintheorem\global\corolcount=0
\gdef\theoremnum{Main Theorem}\else%
\ifx\one\xcorol\ 
\alltheoremnums 
\else\ \alltheoremnums\fi\fi%
\ifthtitle\ \global\thtitlefalse{\rm(\thethtitle)}\fi.%
\hskip1em\bgroup\let\text\rmtext\makeparensRM\it\ignorespaces}

\def\nonumproclaim#1{\vskip-\lastskip
\def\one{#1}\ifx\one\xtheorem\global\corolcount=0\fi
\ifsection\global\sectionfalse\vskip-6pt\fi
\medskip
{\elevensc#1}.\ifx\one\xmaintheorem\global\corolcount=0
\gdef\theoremnum{Main Theorem}\fi\hskip.5pc%
\bgroup\it\makeparensRM\ignorespaces}

\def\endproclaim{\egroup\medskip}


\def\xproof{Proof}
\def\xremark{Remark}
\def\xcase{Case}
\def\xsubcase{Subcase}
\def\xconjecture{Conjecture}
\def\xstep{Step}
\def\xof{of}

\def\deconstruct#1 #2 #3 #4 #5 @{\def\one{#1}\def\two{#2}\def\three{#3}%
\def\four{#4}%
\ifx\two\empty #1\else%
\ifx\one\xproof%
\ifx\two\xof%
  \ifx\three\xcorol Proof of Corollary \rm#4\else%
     \ifx\three\xtheorem Proof of Theorem \rm#4\else\xone\fi%
  \fi\fi%
\else\xone\fi\fi.}

\def\pickup#1 {\def\this{#1}%
\ifx\this\xproof\global\let\go\demoproof
\global\let\enddemo\endproof\else
\ifx\this\xremark\global\let\go\demoremark\else
\ifx\this\xcase\global\let\go\demostep\else
\ifx\this\xsubcase\global\let\go\demostep\else
\ifx\this\xconjecture\global\let\go\demostep\else
\ifx\this\xstep\global\let\go\demostep\else
\global\let\go\demoproof\fi\fi\fi\fi\fi\fi}

\newif\ifnonum
\def\demo#1{\vskip-\lastskip
\ifsection\global\sectionfalse\vskip-6pt\fi
\def\one{#1 }\def\two{#1*}%
\setbox0=\hbox{\expandafter\pickup\one}\expandafter\go\two}

\def\numbereddemo#1{\vskip-\lastskip
\ifsection\global\sectionfalse\vskip-6pt\fi
\def\two{#1*}%
\expandafter\demoremark\two}

\def\demoproof#1*{\medskip\def\xone{#1}
{\ignorespaces\it\expandafter\deconstruct\xone {} {} {} {} {} @%
\unskip\hskip6pt}\rm\ignorespaces}

\def\demoremark#1*{\medskip
{\it\ignorespaces#1\/} \ifnonum\global\nonumtrue\else
 \alltheoremnums\unskip.\fi\hskip1pc\rm\ignorespaces}

\def\demostep#1 #2*{\vskip4pt
{\it\ignorespaces#1\/} #2.\hskip1pc\rm\ignorespaces}

\def\enddemo{\medskip}

\def\endproof{\ifmathqed\global\mathqedfalse\medskip\else
\parfillskip=0pt~~\hfill\qed\medskip
\fi\global\parfillskip0pt plus 1fil\relax
\gdef\enddemo{\medskip}}

\def\qed{\vbox{\hrule\hbox{\vrule height6pt\hskip6pt\vrule}\hrule}}


\def\proofbox{\parfillskip=0pt~~\hfill\qed\vskip1sp\parfillskip=
0pt plus 1fil\relax}







\def\stripbs#1#2*{\def\one{#2}}

\def\emptyspace{ }
\def\nextthing{}
\def\newline{***}
\def\eatone#1{ }

\def\lookatspace#1{\ifcat\noexpand#1\ \else%
\gdef\nextthing{}\xdef\next{#1}%
\ifx\next\emptyspace%
\let\nextthing\emptyspace\else\ifx\next\newline%
\gdef\nextthing{\eatone}\fi\fi\fi\egroup\nextthing#1}

{\catcode`\^^M=\active%
\gdef\spacer{\bgroup\catcode`\^^M=\active%
\let^^M=\newline\obeyspaces\lookatspace}}

\def\ref#1{\seeifdefined{#1}\expandafter\csname\one\endcsname\spacer}

\def\cite#1{\expandafter\ifx\csname#1croref\endcsname\relax[??]\else
\csname#1croref\endcsname\fi\spacer}


\def\seeifdefined#1{\expandafter\stripbs\string#1croref*%
\crorefdefining{#1}}

\newif\ifcromessage
\global\cromessagetrue

\def\crorefdefining#1{\ifdefined{\one}{}
{\ifcromessage\global\cromessagefalse%
\message{\spaces\spaces\spaces\spaces\spaces\spaces\spaces}%
\message{<Undefined reference.}%
\message{Please TeX file once more to have accurate cross-references.>}%
\message{\spaces\spaces\spaces\spaces\spaces\spaces\spaces}\fi[??]}}

\def\label#1#2*{\gdef\ctest{#2}%
\xdef\currlabel{\string#1croref}
\expandafter\seeifdefined{#1}%
\ifx\empty\ctest%
\xdef\labelnow{\write\auxfile{\noexpand\def\currlabel{\the\pageno}}}%
\else\xdef\labelnow{\write\auxfile{\noexpand\def\currlabel{#2}}}\fi%
\labelnow}

\def\ifdefined#1#2#3{\expandafter\ifx\csname#1\endcsname\relax%
#3\else#2\fi}




\def\articlecontents{
\vskip20pt\centerline{\bf Table of Contents}\everypar={}\vskip6pt
\bgroup \leftskip=3pc \parindent=-2pc 
\def\item##1{\vskip1sp\indent\hbox to2pc{##1.\hfill}}}

\def\endcontents{\vskip1sp\leftskip=0pt\egroup}

\def\journalcontents{\vfill\eject
\def\currannalsline{\hfill}
\global\titletrue
\vglue3.5pc
\centerline{\tensc\hskip12pt TABLE OF CONTENTS}\everypar={}\vskip30pt
\bgroup \leftskip=34pt \rightskip=-12pt \parindent=-22pt 
  \def\\ {\vskip1sp\noindent}
\def\pagenum##1{\unskip\parfillskip=0pt\dotfill##1\vskip1sp
\parfillskip=0pt plus 1fil\relax}
\def\name##1{{\tensc##1}}}


\institution{}
\onpages{0}{0}
\def\lastpage{???}
\def\thetitle{Title ???}
\def\theauthors{Authors ???}
\def\thereceived{}
\def\therevised{}

\gdef\split{\relaxnext@\ifinany@\let\next\insplit@\else
 \ifmmode\ifinner\def\next{\onlydmatherr@\split}\else
 \let\next\outsplit@\fi\else
 \def\next{\onlydmatherr@\split}\fi\fi\let\eqnu\xspliteqnu\next}

\gdef\align{\relaxnext@\ifingather@\let\next\galign@\else
 \ifmmode\ifinner\def\next{\onlydmatherr@\align}\else
 \let\next\align@\fi\else
 \def\next{\onlydmatherr@\align}\fi\fi\let\eqnu\xspliteqnu\next}

\def\spliteqnu{{\tenrm\sectandeqnum}\relax}

\def\xspliteqnu{\tag\spliteqnu}

\catcode`@=12

\document

\annalsline{March}{1997}
\startingpage{1}     

\comment
\nopagenumbers
\headline{\ifnum\pageno=1\hfil\else \rightheadline\fi}
\def\rightheadline{\hfil\eightit 
The Macdonald conjecture
\quad\eightrm\folio}

\voffset=2\baselineskip
\endcomment


%
%
%
%
%

\def\iif{\quad\hbox{ if }\quad}

\def\for{\  \hbox{ for } \ }
\def\if{ \ \hbox{ if } \ }

\def\where{\  \hbox{ where } \ }
\def\and{\  \hbox{ and } \ }

\def\equal{\buildrel  def \over =}

\def\la{\lambda}
\def\La{\Lambda}
\def\om{\omega}

\def\th{\theta}
\def\al{\alpha}

\def\ga{\gamma}
\def\ep{\epsilon}

\def\de{\delta}
\def\De{\Delta}
\def\ka{\kappa}

\def\Ga{\Gamma}
\def\ze{\zeta}


\def\tal{\tilde{\alpha}}

\def\tga{\tilde{\gamma}}

\def\tw{\tilde w}

\def\tz{\tilde z}
\def\tb{\tilde b}

\def\hT{\hat{T}}

\def\hw{\hat{w}}

\def\hv{\hat{v}}

\def\C{\bold{C}}

\def\R{\bold{R}}

\def\Z{\bold{Z}}

\def\one{\bold{1}}

\def\0{\bold{0}}

\def\C{\hbox{\bf C}}


\def\l{\Cal{L}}

\def\p{\Cal{P}}

\def\y{\Cal{Y}}
\def\e{\Cal{E}}

\def\x{\Cal{X}}
\def\g{\Cal{G}}

\font\germ=eufb10 
\def\goth#1{\hbox{\germ #1}}

\def\TT{\goth{T}}

\def\AA{\goth{A}}

\font\smm=msbm10 at 12pt 
\def\symbol#1{\hbox{\smm #1}}
\def\lsmash{{\symbol n}}



\title
{Difference Macdonald-Mehta conjecture}

\shorttitle{Macdonald-Mehta conjecture}

\acknowledgements{
Partially supported by NSF grant DMS--9622829}

\author{ Ivan Cherednik}

\institutions{
Math. Dept, University of North Carolina at Chapel Hill,   
 N.C. 27599-3250
\\ Internet: chered\@math.unc.edu
}


\intro 
%
%
%
%
%
\vfil

We formulate and check a difference
counterpart of the Macdonald-Mehta conjecture
and its generalization for the Macdonald polynomials.
It solves the last open problem from the fundamental
paper [M1] and, moreover, gives the formulas for
the Fourier transforms of the
polynomials multiplied by the Gaussian. We also introduce
the reproducing kernel of the difference Fourier transform,
which is an important step towards the difference Harish-Chandra theory.

Mehta suggested a formula for 
the integral of the $\prod_{1\le i< j\le n} (x_i-x_j)^{2k}$
with respect to the Gaussian measure. Macdonald extended it
from $A_{n-1}$ to other root systems and verified his conjecture
for classical ones
by means of Selberg's integrals [M1]. 
It was established by Opdam in [O1] in full generality
using the shift operators. 

The integral is an important normalization constant
for a $k$-deformation of the Hankel transform introduced
by Dunkl [D]. 
The generalized Bessel functions [O3]
multiplied by the Gaussian are
eigenfunctions of this transform. The eigenvalues are given
in terms of this constant. See [D,J] for detail.
The Hankel transform is a
rational degeneration of the Fourier transform in
the Harish-Chandra theory of
spherical functions when the symmetric space $G/K$
is replaced by its tangent space $T_e(G/K)$ with the adjoint
action of $G$
(see [H]).

The harmonic analysis for $G/K$
is much more complicated than that
in the rational case. The 
reproducing kernel of the Fourier transform
is not symmetric, the Gaussian is not Fourier-invariant, and so on.
The zonal spherical functions for dominant weights 
are (trigonometric) polynomials and have no counterparts 
in the rational theory. They play a great role in mathematics
and physics.
When $k=1$ they are
the characters
of finite dimensional representations of $G$.
Unfortunately the  Fourier transform is not very helpful
for the spherical polynomials, but for the characters.

It is the same for the
$k$-deformations, called the
Jack - Heck\-man-Op\-dam polynomials. They  have remarkable
combinatorial properties (Macdonald, Stanley, Hanlon and 
others) and many applications.
A variant of the Mehta-Macdonald 
conjecture in the trigonometric differential setup
is the celebrated Macdonald constant 
term conjecture [M1].
It was proved by Opdam for all root systems
[O1]. 

A difference generalization of the Harish-Chandra
theory was started in [C3,C4]. It fuses together the Fourier
transform and the Gaussian. 
From this viewpoint, it is similar to
the rational case. Moreover, the Fourier transform acts very well on 
the Macdonald $q,t$-polynomials [M2,M3,M4,C4], generalizing
the spherical polynomials, which has no analogue in the 
differential theory. The Macdonald polynomials 
form a basis in the
smallest spherical irreducible representation of the  
double affine Hecke algebra. All spherical representations were
classified in [C1]. They are expected to have 
promising applications in harmonic analysis and combinatorics.

The Macdonald constant term conjecture, as well as the
norm, evaluation, and
duality conjectures, 
were justified in the $q,t$-case 
in [C2,C3,C4]. 
In this paper we complete the theory calculating 
the Fourier transforms of the Macdonald
polynomials multiplied by the Gaussian.
The Gaussian and its transform are proportional.
The formula for the coefficent of proportionality, a difference
counterpart of the Mehta integral,  resembles that for
Macdonald's constant term. It is not surprising because both
are established using similar methods. However
the Gaussian measure is very different from
that due to Macdonald.

In the paper we mainly  consider polynomials and the
pairing based on the constant term. The  Jackson integrals
appear in the last section (see also [C1]).
All statements remain valid for the Jackson
(discrete) pairing. As an application, we
check that the kernel of the Fourier transform
reproduces the Macdonald polynomials, which makes  
its definition self-consistent.

Actually different concepts of integration
do not affect the main formulas 
up to minor renormalizations.
The shift operators always work well.
This holds even when $q$ is a root of unity [C2,C3]
or $k$ is special negative (see [C1],
[DS]). In these cases the Jackson integrals become finite sums,
but it does not change the formulas too much (as long as they
are meaningful).

The author thanks D. Kazhdan and E. Opdam for useful discussion.
The work was partly completed at the University Paris 7.
I am grateful for the kind invitation. 
I acknowledge
my special indebtedness to M. Duflo and P. Gerardin.

%
%
%
\section{Main results}
Let $R=\{\al\}   \subset \R^n$ be a root system of type $A,B,...,F,G$
with respect to a euclidean form $(z,z')$ on $\R^n \ni z,z'$,
$W$ the Weyl group  generated by the the reflections $s_\al$.
We assume that $(\al,\al)=2$ for long $\al$.
Let us  fix the set $R_{+}$ of positive  roots ($R_-=-R_+$), 
the corresponding simple 
roots $\al_1,...,\al_n$, and  their dual counterparts 
$a_1 ,..., a_n,  a_i =\al_i^\vee, \where \al^\vee =2\al/(\al,\al)$.  
The dual fundamental weights
$b_1,...,b_n$  are determined from the relations  $ (b_i,\al_j)= 
\de_i^j $ for the 
Kronecker delta. We will also introduce the dual root system
$R^\vee =\{\al^\vee, \al\in R\}, R^\vee_+$,  the lattices
$$
\eqalignno{
& A=\oplus^n_{i=1}\Z a_i \subset B=\oplus^n_{i=1}\Z b_i, 
}
$$
and  $A_\pm, B_\pm$  for $\Z_{\pm}=\{m\in\Z, \pm m\ge 0\}$
instead of $\Z$. In the standard notations, $A= Q^\vee,\ 
B = P^\vee $ (see [B]).  

Later on, $(\th,\th)=2$ for the maximal root $\th$,
$$
\eqalign{
&\nu_{\al}\ =\ (\al,\al),\  \nu_i\ =\ \nu_{\al_i}, \ 
\nu_R\ = \{\nu_{\al}, \al\in R\}, \cr 
&\rho_\nu\ =\ (1/2)\sum_{\nu_{\al}=\nu} \al \ =
\ (\nu/2)\sum_{\nu_i=\nu}  b_i, \for\al\in R_+,\cr
&r_\nu\ =\ \rho_\nu^\vee \ =\ (2/\nu)\rho_\nu\ =\ 
\sum_{\nu_i=\nu}  b_i,\quad 2/\nu=1,2,3. 
} 
\eqnu
$$

We will mainly use $r_\nu$ and $r=\sum_{\nu}r_\nu$ in the paper.
The theory depends on the parameters
$q, t_\nu, \nu\in \nu_R$. It is convenient to set
$$
q_\nu=q^{2/\nu},\ t_\nu=q_\nu^{k_\nu},\
q_\al=q_{\nu}, t_\al=t_\nu \for \nu=\nu_\al \and
r_k=\sum_\nu k_\nu r_\nu.
$$

Let us formally put 
$$
\eqalign{
&x_i=q^{b_i},\  x_b=q^b= \prod_{i=1}^n
x_i^{l_i} \for b=\sum_{i=1}^n l_i b_i,
} 
\eqnu
$$
and introduce the algebra
$\C(q,t)[x]$  of polynomials in terms of $x_i^{\pm 1}$ with the
coefficients belonging to the field $\C(q,t)$ of rational functions.

The coefficient of $x^0=1$ ({\it the constant term})
will be denoted by $\langle \  \rangle$. The following product
(the Macdonald truncated $\theta$-function) is a 
Laurent series in $x$ with coefficients in  the algebra
$\C[t][[q]]$ of formal (holomorphic) series in $q$ over polynomials in $t$:
$$
\eqalign{
&\De\ =\ \prod_{\al \in R_+}
\prod_{i=0}^\infty {(1-x_aq_\al^{i}) (1-x_a^{-1}q_\al^{i})
\over
(1-x_a t_\al q_a^{i}) (1-x_a^{-1}t_\al^{}q_\al^{i})},
\ a=\al^\vee.
}
\eqnu
\label\Demu\eqnum*
$$
Here 
$(1-(\cdot))^{-1}$ are replaced by $1+(\cdot)+(\cdot)^2+\ldots\ .$
We note that  $\De\in 
\C(q,t)[x]$ if $t_\nu=q_\nu^{k_\nu}$ for $k_\nu\in \Z_+$.

By the {\it Gaussians} $\tga^{\pm 1}$ we mean 
$$
\eqalign{
&\tga= \sum_{b\in B} q^{-(b,b)/2}x_b,\ 
\tga^{-1}= \sum_{b\in B} q^{(b,b)/2}x_b.
}
\eqnu
\label\gauser\eqnum*
$$
The multiplication by $\tga^{-1}$ preserves 
the space of Laurent series
 with coefficients from $\C[t][[q]]$.

\proclaim{Macdonald-Mehta Theorem}
$$
\eqalignno{
&\langle \tga^{-1}\De\rangle\ =\ 
|W|\prod_{\al\in R_+}\prod_{ j=0}^{\infty}\Bigl({ 
1- q_\al^{(r_k,\al)+j}\over
1-t_\al q_\al^{(r_k,\al)+j} }\Bigr).
&\eqnu
\label\mehta\eqnum*
}
$$
\label\MEHTA\theoremnum*
\endproclaim

Here $t$ is either a formal parameter or
the right hand side is to be  understood as the
corresponding limit when some $k_\nu \in \Z_-$. 

The {\it monomial symmetric functions} 
$m_{b}\ =\ \sum_{c\in W(b)}x_{c}$ for $b\in B_-$
form a basis of the space 
 $\C[x]^W$ of all $W$-invariant polynomials. 
We  introduce the {\it Macdonald 
polynomials} $p_b(x),\   b \in B_-$, by means of 
the conditions
$$
\eqalignno{
&p_b-m_b\ \in\ \oplus_c\C(q,t)m_{c},\
\langle p_b m_{c}\De\rangle = 0 \for c> b,  &\eqnu 
 \cr
&\hbox{where\ }
 c\in B_-,\ 
 c> b \hbox{\ means\ that\ }  c-b \in A_+, c\neq b.
\label\macd\eqnum*
}
$$
They can be determined by the Gram - Schmidt process 
(see [M2,M3])  
 and form a 
basis in $\C(q,t)[x]^W$. As it was established by Macdonald,
they are pairwise orthogonal
for the pairing $\langle f(x)g(x^{-1})\De\rangle$.
We consider $t$ as a formal parameter. If it is a number then it is
necessary to avoid certain negative $k$.

\proclaim{ Theorem}
Given $b,c\in B_-$ and the corresponding Macdonald
polynomials $p_b, p_c$,
$$
\eqalign{
&\langle
p_b p_c \tga^{-1}\De
\rangle = 
 q^{(b,b)/2+(c,c)/2 -(b+c,r_k)} 
p_c(q^{b-r_k})p_b(q^{r_k})\langle \tga^{-1}\De\rangle \cr
&= q^{(b,b)/2+(c,c)/2 -(c,r_k)} p_c(q^{b-r_k})
\, |W|\prod_{\al\in R_+}\, \prod_{ j=-(\al,b)}^{\infty}\Bigl({ 
1- q_\al^{(r_k,\al)+j}\over
1-t_\al q_\al^{(r_k,\al)+j}}\Bigr),
}
\eqnu
\label\planch\eqnum*
$$
where $x_c(q^b)\equal q^{(b,c)}$.
\label\PLANCH\theoremnum*
\endproclaim           

In the last formula, we
used the Macdonald evaluation conjecture
proved in [C3]:
$$
\eqalignno{
&p_b(q^{r_k})\ =\ 
q^{(r_k,b)}\prod_{\al\in R_+}
\prod_{ j=1}^{-(\al,b)}  
\Bigl(
{ 
1- t_\al q_\al^{(r_k,\al)+j-1}
 \over
1-  q_\al^{(r_k,\al)+j-1}
}
\Bigr). &\eqnu
\label\value\eqnum*
}
$$
To make the products meaningful we expand the
coefficients of the polynomials $p_b$ in terms of
$q$. They are from $\C[t][[q]]$.

The formula is an important
particular case of the difference Fourier-Plancherel theorem. 
We will extend it to 
the non-symmetric Macdonald polynomials.
More general results will be considered in the next paper.

There is a straightforward passage
to non-reduced root systems.
One may also take $x_\al$ in place of $x_a$,
changing $q_\al$ by $q$, choosing
$c$ from the lattice $P$, and substituting
$(r_k,\al)\to(\rho_k,\al^\vee)$,
$(b,r_k)\to (b,\rho_k)$ 
in (\ref\Demu), (\ref\mehta) and (\ref\planch), (\ref\value).
The Gaussian remains the same (for $B$).

A rational-differential counterpart of (\ref\planch)
was verified by Dunkl and de Jeu
(see [D], Theorem 3.2). 
Theorem \ref\PLANCH
is the cornerstone of the 
theory of the Fourier transform (cf. [J], Lemma 4.11).
Generally speaking,
the latter is the map $f(x)\to \hat{f}(\la)=\int p_\la(x)f(x)\De$,
where $p_\la$ is a $W$-invariant eigenfunction of the generalized
Macdonald operators for a proper choice of the 
integration and $\la$. For $\int=\langle \rangle$,
(\ref\planch) determines the
Fourier transform and its inverse in the space of
$W$-symmetric Laurent polynomials
multiplied by the Gaussian. This space is identified with a subspace
of all functions in 
$\la= q^{b-r_k}, b\in B$.

Both statements are new.
In the case of rank one, (\ref\mehta) resembles
the so-called quintuple product identity and 
the formulas from [AW] (the $BC_1$ case).
So it is likely to be related to the known one-dimensional identities.

The definition of the {\it Jackson integral} for
$W$-symmetric $f$ is as follows: $\langle f\rangle_\xi =
\sum_{a\in B} f(q^{\xi+a})$. The integrals 
are formal series in terms of $q$ or
 functions of $q,t,\xi$ when  $|q|<1$.   
We choose $\ga(q^z)=q^{(z,z)/2}$ and 
$$
\eqalign{
&\De^\circ\  =\  \prod_{\al \in R_+}
\prod_{i=1}^\infty {(1-x_a t_\al^{-1}q_\al^{i}) 
(1-x_a^{-1}t_\al^{-1}q_\al^{i})
\over
(1-x_a q_a^{i}) (1-x_a^{-1}q_\al^{i})}.
}
\eqnu
\label\Decirc\eqnum*
$$
For instance, $\langle \ga\rangle_\xi= \sum_{a\in B} q^{(\xi+a,\xi+a)/2}=
\tga^{-1}(q^\xi)q^{(\xi,\xi)/2}$. We assume 
that $\De^\circ(q^{\xi+b})$ is
well-defined, so $(\al,\xi)\not\in \Z$ for all $\al\in R_+$.

\proclaim{ Theorem}
Given $b,c\in B_-$ and the corresponding Macdonald
polynomials $p_b, p_c$,
$$
\eqalignno{
&\langle
p_b(x) p_c(x^{-1})\,\ga\De^{\circ}
\rangle_\xi\  = \cr
& q^{-(b,b)/2-(c,c)/2 +(b+c,r_k)} 
p_c(q^{b-r_k})p_b(q^{r_k})\langle \ga\De^\circ\rangle_\xi,
&\eqnu \cr
\label\planjack\eqnum*
&
\langle \ga\De^\circ\rangle_\xi\ =\ 
\langle \ga\rangle_\xi \prod_{\al\in R_+}\prod_{ j=1}^{\infty}\Bigl({ 
1- t_\al^{-1}q_\al^{-(r_k,\al)+j}\over
1-q_\al^{-(r_k,\al)+j} }\Bigr).
&\eqnu
\label\mehjack\eqnum*
}
$$
\label\PLANJACK\theoremnum*
\endproclaim

In this formulas $t$ is arbitrary provided the existence
of $p_{b,c}$.  The right hand side of (\ref\mehjack) is 
considered as the corresponding limit if $k_\nu\in \Z_+\setminus\{0\}$.

%
%
%
%
\section {Affine Weyl groups} 
The vectors $\ \tal=[\al,k] \in 
\R^n\times \R \subset \R^{n+1}$ 
for $\al \in R, k \in \Z $ 
form the {\it affine root system} 
$R^a \supset R$ ( $z\in \R^n$ are identified with $ [z,0]$).
We add  $\al_0 \equal [-\th,1]$ to the  simple roots 
for the  maximal root $\th \in R$.
The corresponding set $R^a_+$ of positive roots coincides
with $R_+\cup \{[\al,k],\  \al\in R, \  k > 0\}$. 

We denote the Dynkin diagram and its affine completion with
$\{\al_j,0 \le j \le n\}$ as the vertices by $\Ga$ and $\Ga^a$.
The set of
the indices of the images of $\al_0$ by all 
the automorphisms of $\Ga^a$ will be denoted by $O$ ($O=\{0\} 
\for E_8,F_4,G_2$). Let $O^*={r\in O, r\neq 0}$.
The elements $b_r$ for $r\in O^*$ are the so-called minuscule
weights ($(b_r,\al)\le 1$ for
$\al \in R_+$).

Given $\tal=[\al,k]\in R^a,  \ b \in B$, let  
$$
\eqalignno{
&s_{\tal}(\tz)\ =\  \tz-(z,\al^\vee)\tal,\ 
\ b'(\tz)\ =\ [z,\ze-(z,b)]
&\eqnu
}
$$
for $\tz=[z,\ze] \in \R^{n+1}$.

The {\it affine Weyl group} $W^a$ is generated by all $s_{\tal}$
(we write $W^a = <s_{\tal}, \tal\in R_+^a>)$. One can take
the simple reflections $s_j=s_{\al_j}, 0 \le j \le n,$ as its
generators and introduce the corresponding notion of the  
length. This group is
the semi-direct product $W\lsmash A'$ of 
its subgroups $W=<s_\al,
\al \in R_+>$ and $A'=\{a', a\in A\}$, where
$$
\eqalignno{
& a'=\ s_{\al}s_{[\al,1]}=\ s_{[-\al,1]}s_{\al}\for a=\al^{\vee},
\ \al\in R.
&\eqnu
}
$$

The {\it extended Weyl group} $ W^b$ generated by $W\and B'$
(instead of $A'$) is isomorphic to $W\lsmash B'$:
$$
\eqalignno{
&(wb')([z,\ze])\ =\ [w(z),\ze-(z,b)] \for w\in W, b\in B.
&\eqnu
}
$$
Later on  $b$ and $b'$ will not be distinguished. 

Given $ b\in B$, the decomposition $b= \pi_b\om_b,
\om_b \in W$
can be uniquely determined  from the  condition:
$\om_b(b) = b_-\in B_-$ where the length $l(\om_b)$
of $\om$ in terms of $\{s_1,\ldots,s_n\}$  is the smallest possible. 
For instance, let $\pi_r=\pi_{b_r}, r \in O$. 
They leave $\Ga^a$ invariant 
and form a group denoted by $\Pi$, 
which is isomorphic to $B/A$ by the natural 
projection $\{b_r \to \pi_r\}$. As to $\{\om_r\}$,
they preserve the set $\{-\th,\al_i, i>0\}$.
The relations $\pi_r(\al_0)= \al_r= (\om_r)^{-1}(-\th)$ distinguish the
indices $r \in O^*$. Moreover (see e.g. [C2]):
$$
\eqalignno{
& W^b  = \Pi \lsmash W^a, \where
  \pi_rs_i\pi_r^{-1}  =  s_j \if \pi_r(\al_i)=\al_j,\  0\le j\le n.
&\eqnu
}
$$

We extend   
the  length to $W^b$.
Given $r\in O^*,\  \tw \in W^a$, and a reduced  
decomposition $\tw\ =\ s_{j_l}...s_{j_2} s_{j_1} $ with respect to
$\{s_j, 0\le j\le n\}$, we call $l\ =\ l(\hw)$ the {\it length} of 
$\hw = \pi_r\tw \in W^b$. Similarly, $l_\nu(\hw)$ is the number
of $s_j$ with $\nu_j=\nu$.

Let us introduce a partial ordering
on $B$. Here and further $b_{-}$ is the unique elements
from $B_{-}$ which belong to the orbit $W(b)$. Namely,
$b_-=\om_b(b)$.
So the equality   $c_-=b_- $ means that $b,c$
belong to the same orbit.
Set 
$$
\eqalignno{
&b \le c, c\ge b \for b, c\in B \iif c-b \in A_+,
&\eqnu 
\label\order\eqnum*
\cr
&b \preceq c, c\succeq b \iif b_-< c_- \hbox{\ \ or\ \ }
b_-=c_- \hbox{\ and\ } b\le c.
&\eqnu
\label\succ\eqnum*
}
$$
We  use $<,>,\prec, \succ$ respectively if $b \neq c$.

%
%
%
\section{ Difference operators}
We put    
$m=2 \for D_{2k} \and C_{2k+1},\ m=1 \for C_{2k}, B_{k}$,
otherwise $m=|\Pi|$. We use 
the parameters $ q, t_\nu=q_\nu^{k_\nu}, q_\nu=q^{2/\nu} (\nu \in \nu_R)$
 and the variables $x_1,\ldots,x_n$ are from Section 1,
$$
\eqalignno{
&   t_{[\al,k]} = t_{\al}=t_{\nu_\al},\ t_j = t_{\al_j},
\where [\al,k] \in R^a,\ 0\le j\le n, \cr 
& x_{\tb}\ =\ \prod_{i=1}^nx_i^{l_i} q^{ k} 
\if \tb=[b,k],
&\eqnu \cr 
&\hbox{for\ } b=\sum_{i=1}^n l_i b_i\in B,\ k \in {1\over m}\Z.
}
$$
\label\Xde\eqnum*

We will also consider polynomials in $q^{\pm 1/m}$ and  $\{t_\nu^{\pm 1/2} \}$,
using the notation $\C[q^{\pm 1/m},t^{\pm 1/2}]$. 
The elements $\hw \in W^b$ act in $\C[q^{\pm 1/m}][x]$
 by the 
formulas:
$$
\eqalignno{
&\hw(x_{\tb})\ =\ x_{\hw(\tb)}. 
&\eqnu}
$$
 In particular:
$$
\eqalignno{
&\pi_r(x_{b})\ =\ x_{\om^{-1}_r(b)} q^{(b_{r^*},b)} 
\for \al_{r^*}\ =\ \pi_r^{-1}(\al_0), \ r\in O^*.
&\eqnu}
$$
\label\pi\eqnum*

The {\it Demazure-Lusztig operators} (see  
[KL, KK, C2]) 
$$
\eqalignno{
&T_j\  = \  t_j ^{1/2} s_j\ +\ 
(t_j^{1/2}-t_j^{-1/2})(x_{a_j}-1)^{-1}(s_j-1),
\ 0\le j\le n.
&\eqnu
\label\Demaz\eqnum*
}
$$
preserve $\C[q^{\pm 1/m},t^{\pm 1/2}][x]$.
We note that only $T_0$ involves $q$: 
$$
\eqalign{
&T_0\  =  t_0^{1/2}s_0\ +\ (t_0^{1/2}-t_0^{-1/2})
( q X_{\th}^{-1} -1)^{-1}(s_0-1),\cr
&\where
s_0(X_i)\ =\ X_iX_{\th}^{-(b_i,\th)} q^{(b_i,\th)}.
}
\eqnu
$$

Given $\tw \in W^a, r\in O,\ $ the product
$$
\eqalignno{
&T_{\pi_r\tw}\equal \pi_r\prod_{k=1}^l T_{i_k},\where 
\tw=\prod_{k=1}^l s_{i_k},
l=l(\tw),
&\eqnu
\label\Tw\eqnum*
}
$$
does not depend on the choice of the reduced decomposition
of $\tw$
(because $\{T\}$ satisfy the same `braid' relations as $\{s\}$ do).
Moreover,
$$
\eqalignno{
&T_{\hv}T_{\hw}\ =\ T_{\hv\hw}\  \hbox{ whenever}\ 
 l(\hv\hw)=l(\hv)+l(\hw) \for
\hv,\hw \in W^b.
&\eqnu
\label\TT\eqnum*
}
$$

  In particular, we arrive at the pairwise 
commutative operators (see [C2]):
$$
\eqalignno{
& Y_{b}\ =\  \prod_{i=1}^nY_i^{k_i} \if  
b=\sum_{i=1}^nk_ib_i\in B,\where  
 Y_i\equal T_{b_i},
&\eqnu
\label\Yb\eqnum*
}
$$
satisfying the relations
$$
\eqalign{
&T^{-1}_iY_b T^{-1}_i\ =\ Y_b Y_{a_i}^{-1} \if (b,\al_i)=1,
\cr
& T_iY_b\ =\ Y_b T_i \if (b,\al_i)=0, \ 1 \le i\le  n.
}
\eqnu
$$

%
%
%
%
\section{ Macdonald polynomials}
Recall that $\langle  f \rangle$ is the constant term
of $f$. We will switch from $\De$ to 
$$
\eqalign{
&\mu\ = \mu^{(k)}=\prod_{\al \in R_+}
\prod_{i=0}^\infty {(1-x_aq_\al^{i}) (1-x_a^{-1}q_\al^{i+1})
\over
(1-x_a t_\al q_a^{i}) (1-x_a^{-1}t_\al^{}q_\al^{i+1})},\
a=\al^\vee.
}
\eqnu
\label\mu\eqnum*
$$
It is considered as a Laurent series with the
coefficients in  $\C[t][[q]]$.

Let
$\mu_1\equal \mu/\langle \mu \rangle$, 
where the formula for the
constant term of $\mu$ is as follows
(see [C2]):
$$
\eqalign{
&\langle\mu\rangle\ =\ \prod_{\al \in R_+}
\prod_{i=1}^{\infty} { (1- q_\al^{(r_k,\al)+i})^2
\over
(1-t_\al q_\al^{(r_k,\al)+i}) (1-t_\al^{-1}q_\al^{(r_k,\al)+i})
}.
}
\eqnu
\label\consterm\eqnum*
$$
It is a Laurent series with coefficients 
in $\C(q,t)$, and 
$\mu_1^*\ =\ \mu_1$  with respect to the involution 
$$
 x_b^*\ =\  x_{-b},\ t^*\ =\ t^{-1},\ q^*\ =\ q^{-1}.
$$
 
Setting 
$$
\eqalignno{
&\langle f,g\rangle_1\ =\langle \mu_1 f\ {g}^*\rangle_1\ =\ 
\langle g,f\rangle_1^* \for
f,g \in \C(q,t)[x], 
&\eqnu
\label\innerpro\eqnum*  
}
$$
we  introduce the {\it non-symmetric Macdonald 
polynomials} $e_b(x)=e_b^{(k)},\   b \in B$, by means of 
the conditions
$$
\eqalignno{
&e_b-x_b\ \in\ \oplus_{c\succ b}\C x_c,\
\langle e_b, x_{c}\rangle_1 = 0 \for B_-\ni c\succ b.
&\eqnu
\label\macd\eqnum*
}
$$
They are well-defined
because the  pairing   
is non-degenerate 
 and form a 
basis in $\C(q,t)[x]$. 

Replacing here $\mu_1$ by $\mu$ and involving the norm-formulas
from [C4] (the Main Theorem), we establish that 
the expansions of the coefficients of $e_b$ in terms of $q$ belong
to $\C[t^{\pm 1}][[q]]$. Applying $^*$, these coefficients also belong to
 $\C[t^{\pm 1}][[q^{-1}]]$ when expanded in $q^{-1}$. The same holds
for $e_b^*$.

This definition is due to Macdonald (for
 $t_\nu=q_\nu^k,\ k\in \Z_+$),
 who extended 
 Opdam's non-symmetric polynomials introduced
in the degenerate (differential) case in [O2].
The general case was considered in [C4]. 
Another approach is based on the $Y$-operators
(see [M4],[C4]):

\proclaim {Proposition}
 The polynomials $\{e_b, b\in B\}$
 are eigenvectors of
 the operators $\{L_f\equal f(Y_1,\cdots, Y_n), f\in \C[x]\}$:
$$
\eqalignno{
&L_{f}(e_b)\ =\ f(q^{-b_\#})e_b, \where
b_\#\equal b-\om_b^{-1}(r_k),
&\eqnu
\label\Yone\eqnum*
\cr
& x_a(q^{b_\#})\ =\
q^{(a,b)}\prod_\nu t_\nu^{-(\om_b^{-1}(\rho_\nu),a)},
\om_b \hbox{\ is\ from\ Sec.\ 2.}
&\eqnu
\label\xaonb\eqnum*
}
$$
\endproclaim
\label\YONE\theoremnum*

\comment
Actually here we need a 
 bigger  space.
Namely, in the space $\l_q$ of  Laurent series $f=\sum_{b\in B} c_b x_b$
with the coefficients $c_b\in \C[t][[q]]$ such
that $0\le m_b\le C_f|(b,r)|$ for a constant $C_f> 0$
depending on $f$.
This space is invariant with respect to
the action of the  affine Weyl group. We may also
multiply its elements by the Laurent polynomials from 
$\C(q,t)[x]$ and by $\tga^{-1}$.
\endcomment

Similarly, the symmetric Macdonald polynomials $p_b=p_b^{(k)}$ from 
(\ref\macd) are eigenfunctions of the $W$-invariant
operators $L_f=f(Y_1,\cdots,Y_n)$ for $f\in \C(q,t)[x]^W$:
$$
\eqalignno{
&L_{f}(p_b)\ =\ f(q^{-b+r_k} ) p_b,\
b\in B_-.
&\eqnu
\label\Lf\eqnum*
}
$$

They are connected with $\{e\}$ as follows (see [M4,C4] and [O2] in the
differential case):
$$
\eqalign{
&p_{b}\ =\ \p_b^t e_{b},\ \ b=b_-\in B_-,\cr 
&\p_b^t\equal\sum_{c\in W(b)}
\prod_\nu t_\nu^{l_\nu(w_c)/2} \hT_{w_c},
}
\eqnu
\label\symmetr\eqnum*
$$
where $w_c\equal \om_c^{-1}w_0$ for the longest $w_0\in W$.

Following [C2-C3], let us 
fix a subset $v\in \nu_R$ and introduce the 
{\it shift operator} by the formula:
$$
\eqalignno{
&\g_v \ =\ \g_v^{(k)}\ =\ 
(\x_v) ^{-1}\y_v,
 &\eqnu
\label\shift\eqnum*
}
$$
$$
 \x_v  = \prod_{\nu_\al\in v}((t_\al x_{a})^{1/2}-
(t_\al x_{a})^{-1/2}),\  \y_v  = \prod_{\nu_\al\in v}
(t_\al Y_{a}^{-1})^{1/2}-
(t_\al Y_{a}^{-1})^{-1/2}).
$$
Here $a=\al^\vee, \al\in R_+,$
$\x_{v}= \x_{v}^{(k)}$ and $\y_{v}=\y_{v}^{(k)}$ 
belong to $\C[t^{\pm 1/2}] [x]$ and
$\C[t^{\pm 1/2}] [Y]$ respectively.

\proclaim{Proposition}
The operators ${\g}_{v}$
 are $W$-inva\-riant and pre\-serve the space 
$\C[q^{1/m},t^{\pm 1/2}][x]^W$. 
If $t_{\nu}=1$ for 
$\nu\not\in v$ then
$$
\eqalign{
&\g_v^{(k)} (p_{b}^{(k)})= g_v^{(k)}(b)
p_{b+r_v}^{k+v} \for\cr
&g_v^{(k)}(b)\ =\ 
\prod_{\al\in R_+,\nu_\al\in v} (q_\al^{(r_k-b,\al)/2} -
t_\al q_\al^{(b-r_k,\al)/2}),
}
\eqnu
\label\Gnu\eqnum*
$$
where $r_v=\sum_{\nu\in v}r_\nu,\ k+v=\{k_\nu+1 , k_{\nu'}\}$
for $\nu\in v\not\ni \nu'\ $, $p_c=0 \for c\not\in B_-$.
\label\GP\theoremnum*
\endproclaim

%
%
%
%
\section { Fourier transforms} 
Proofs of the following theorems are based on the
analysis of the automorphisms of the double affine Hecke
algebras. The technique is similar to
that from  [C1-C4] and will be exposed in more detail
in the next paper.

We will mainly use the {\it renormalized Macdonald polynomials} :
$$
\eqalignno{ 
&\ep_{b}=e_b/e_b(q^{-r_k})
= q^{-(r_k,b_-)}\prod_{[\al,j]\in \Lambda_b}
\Bigl(
{ 
1-  q_\al^{(r_k,\al)+j}
 \over
1- t_\al q_\al^{(r_k,\al)+j}
}
\Bigr)\ e_b,
&\eqnu\cr
\label\epsil\eqnum*
&\Lambda_b=\{[\al,j],\ \al\in R_+, 0<j<-(\al,b_-)\if (\al,b)>0,
&\eqnu\cr
\label\jalset\eqnum*
&0<j\le -(\al,b_-)\if (\al,b)<0\},\ \ b_-=\om_b(b).  
}
$$
Here we applied the Main Theorem from [C4].
This normalization is very convenient
in the difference harmonic analysis.
For instance, the duality relations are especially simple:
$\ep_b(q^{c_\#})=\ep_c(q^{b_\#})$.

\proclaim{Theorem}
Given $b,c\in B$ and the corresponding renormalized
polynomials $\ep_b, \ep_c$,
$$
\eqalignno{
&\langle \ep_b \ep_c \tga^{-1}\mu
\rangle\ =\ 
q^{(b_\#,b_\#)/2+(c_\#,c_\#)/2 -(r_k,r_k)} 
\ep_c(q^{b_\#})\langle \tga^{-1}\mu\rangle,
&\eqnu\cr
\label\epep\eqnum*
&\langle
\ep_b \ep_c^* \tga^{-1}\mu
\rangle\ =\ 
q^{(b_\#,b_\#)/2+(c_\#,c_\#)/2 -(r_k,r_k)} 
\ep_c^*(q^{b_\#})\langle \tga^{-1}\mu\rangle,
&\eqnu\cr
\label\epeps\eqnum*
&\langle
\ep_b \ep_c^* \tga\mu^*
\rangle\ =\ 
q^{-(b_\#,b_\#)/2-(c_\#,c_\#)/2 +(r_k,r_k)} 
\ep_c(q^{b_\#})\langle \tga\mu^*\rangle.
&\eqnu
\label\epepl\eqnum*
}
$$
\label\EPEP\theoremnum*
\endproclaim


The first two formulas make sense when we expand 
$\ep_b,\ep_c,\ep_c^*$
in terms of $q$. The coefficients of $\ep,\ep^*$ belong to $\C[t^{\pm 1/2}][[q]]$.
In the last formula, $\mu^*$ is expanded in
terms of powers of $q^{-1}$, as well as  $\ep_b,\ep_c^*$. Thus
(\ref\epepl) results from (\ref\epeps) and the duality.
The product $\mu\tga^{-1}$  generalizes
the (radial) Gaussian measure  in the theory of Lie groups and
symmetric spaces. Actually (\ref\epepl) is a formula for
the {\it Fourier transform} of $\ep_c\tga^{-1}$ (see [C3]). 

The following theorem
is equivalent to Theorem \ref\EPEP  although it does not
involve $c$ and any scalar products.
We will use the conjugation `$\iota$' :  $q^\iota= q^{-1},
\ t^\iota= t^{-1},\ x_b^\iota= x_b$.

\proclaim{Theorem}
Given $b\in B$,
$$
\eqalignno{
&\ep_b(Y_1^{-1},\cdots,Y_n^{-1}) (\tga)\ =\ 
q^{(b_\#,b_\#)/2 -(r_k,r_k)/2} 
\ep_b\tga,
&\eqnu\cr
\label\epy\eqnum*
&\ep_b^\iota(Y_1,\cdots,Y_n) (\tga)\ =\ 
q^{(b_\#,b_\#)/2 -(r_k,r_k)/2} 
\ep_b^*\tga,
&\eqnu\cr
\label\epys\eqnum*
&\ep_b^\iota(Y_1,\cdots,Y_n) (\tga^{-1})\ =\ 
q^{-(b_\#,b_\#)/2 +(r_k,r_k)/2} 
\ep_b\tga^{-1}.
&\eqnu
\label\epyl\eqnum*
}
$$
\label\EPY\theoremnum*
\endproclaim

The first formula holds for any normalizations of $e_b$.
Moreover, since the coefficient of proportionality
$q^{(b_\#,b_\#)/2 -(r_k,r_k)/2}$ is the same for all
$b$ from the same  $W$-orbit, it can be applied to linear combinations
of $e_c, c\in W(b)$. For instance,
will can use it for the symmetric Macdonald polynomials
or for the $t$-antisymmetric ones (the next section).

Actually we do not need a presentation of
 $\tga^{\pm 1}$ as a Laurent series in this theorem.
One may take  $\ga^{\pm 1}\ =\  
q^{\pm\Sigma_{i=1}^n b_i {\al_i}/2}$ or  any
 $W$-invariant solution of the following 
system of difference equations:
$$
\eqalign{
&b_j(\ga)\ =\  q^{(1/2)\Sigma_{i=1}^n (b_i-(b_j,b_i))
({\al_i}- \de_i^j)}\ =\cr
& \ga  q^{-b_j+ (b_j,b_j)/2 }\ =\  q^{(b_j,b_j)/2} x_j^{-1}\ga, \cr 
&b_j(\ga^{-1})\  =\  q^{-(b_j,b_j)/2} x_j\ga^{-1} 
\for 1\le j\le n.
}
\eqnu
\label\gausseq\eqnum*
$$

The formula (\ref\epyl)  for the symmetric polynomials 
was verified in [C3] via the
roots of unity. See (4.19) and the end of the proof
of Corollary 5.4 (use that
$p_b^\iota=p_b$ for the symmetric polynomials). 
The same method works well
in the non-symmetric case and for other
identities.
A more natural proof will appear in the next paper.
To make the formulas complete we need to know the
 coefficent
of proportionality:

\proclaim{Theorem}
$$
\eqalignno{
&\langle \tga^{-1}\mu\rangle\ =\ 
\prod_{\al\in R_+}\prod_{ j=1}^{\infty}\Bigl({ 
1- q_\al^{(r_k,\al)+j}\over
1-t_\al q_\al^{(r_k,\al)+j} }\Bigr).
&\eqnu
\label\mehtamu\eqnum*
}
$$
\label\MEHTAMU\theoremnum*
\endproclaim
When $k_\nu^0\in \Z_-\setminus \{0\}$, by the right hand side we
mean the limit as $k_\nu\to k_\nu^0$.

There is one more reformulation of Theorem \ref\EPEP.
Let us introduce another set of variables
$\la_i, 1\le i\le n,$ and the corresponding $\tga_\la^{\pm 1}$.
We will also use the operators $Y_b^\la, T_i^\la$ acting with respect
to $\la$ by the same formulas as for $x$.
The coefficients of the
following Laurent series in terms of $\{x_i,\la_j\}$
are well-defined:
$$
\eqalignno{
& q^{(r_k,r_k)/2}\e_{q}(x,\la)\tga_x^{-1}\tga_\la^{-1}\ =\
\sum_{b\in B} q^{(b_\#,b_\#)/2 -(r_k,r_k)/2} 
\ {\ep_b^*(x)\ \ep_b(\la)\over
\langle \ep_b,\ep_b\rangle_1},
&\eqnu\cr
\label\exsela\eqnum*
& q^{-(r_k,r_k)/2}\e_{q^{-1}}(x,\la)\tga_x\tga_\la\ =\
\sum_{b\in B} q^{-(b_\#,b_\#)/2 +(r_k,r_k)/2} 
\ {\ep_b(x)\ \ep_b(\la)\over
\langle \ep_b,\ep_b\rangle_1}.
&\eqnu\cr
\label\exela\eqnum*
}
$$
They belong respectively to
$\C[t^{\pm 1/2}][[q^{1/m}]]$ and 
$\C[t^{\pm 1/2}][[q^{-1/m}]]$.
The latter series is obviously symmetric: 
 $\e_{q^{-1}}(x,\la)=\e_{q^{-1}}(\la,x)$.
It holds for the first one too, but requires some proving
(use (\ref\epep) instead of (\ref\epeps)).

\proclaim{Theorem} The series $\e_q(x,\la)$ is symmetric: 
$$
\eqalignno{
& q^{(r_k,r_k)/2}\e_{q}(x,\la)\tga_x^{-1}\tga_\la^{-1}\ =\ 
\sum_{b\in B} q^{(b_\#,b_\#)/2 -(r_k,r_k)/2} 
\ {\ep_b(x)\ \ep_b^*(\la)\over
\langle \ep_b,\ep_b\rangle_1 }.
&\eqnu\cr
\label\exelas\eqnum*
}
$$
Moreover,
$$
\eqalign{
& Y^x_b(\e(x,\la))= (\la_b)^{-1}\e(x,\la),\ 
 Y^\la_b(\e(x,\la))= (x_b)^{-1}\e(x,\la),\cr 
&T_i^x(\e(x.\la))= T_i^\la(\e(x,\la)), \ 1\le i\le n,\  
\for \e=\e_{q^{\pm 1}}.
}
\eqnu
\label\yelet\eqnum*
$$
\label\EGAGA\theoremnum*
\endproclaim

We note that the right hand side of (\ref\exsela)
is a holomorphic function for all $x_i\neq 0\neq \la_j$
and generic $t$
provided that $|q|< 1$. Similarly, (\ref\exela) is holomorphic
for non-zero $x_i,\la_j$ when $|q|> 1$.
It is checked by means of the recurrence relations (the Pieri rules)
from [C4], Theorem 5.4. Dividing the Gaussians out, we see that 
$\e(x,\la)$ are meromorphic with the poles coming from the zeros
of $\tga^{\pm 1}$. The latter set (of divisors) is related
to the Macdonald identities (product formulas).
These functions are also invariant with respect to the intertwiners
from Section 5, [C1]. It readily results from the definition and the
theorem.
The functions $\e(x,\la)$ generalize the reproducing kernel of the
Dunkl transform (which is also symmetric) and that from the 
Harish-Chandra theory (which is not). The passage to the $p$-polynomials
is a straightforward $t$-symmetrization.
\footnote {*}{ 
In a recent paper by Baker, Forrester (q-alg/9701039, 30 Jan 1997),
a certain  generalization of the Dunkl kernel was
suggested (formula (5.1)) in the case of $GL_n$. It is a certain
sum in terms of the non-symmetric polynomials resembling our 
(\ref\exelas). However it is not $x\leftrightarrow \la$
symmetric (according to the authors), 
and there is no counterpart of
the important formula $ Y^\la_b(\e(x,\la))= (x_b)^{-1}\e(x,\la)$ 
in the paper.
The meaning (existence, convergence) of the infinite sums 
is not discussed by the authors. In our approach we 
separate the contribution of the Gaussians. 
After this the analysis is mainly due to [C4].
}

%
%
%
\section{ Mehta integral}
Following [C2,C3] we will use the shift operator to 
verify Theorem \ref\MEHTAMU
in a way similar to that in the differential case [O1].

First we take $k_\nu\in \Z_+$ and replace $\mu$, which
is a trigonometric polynomial in this case,  by a proportional one:
$$
\eqalignno{
\widetilde{\mu}=\widetilde{\mu}^{(k)}\equal \prod_{\al\in R_+}
&((q_\al^{k_\al-1} x_a)^{1/2}-(q_\al^{k_\al-1} x_a)^{-1/2})\cr
\cdots ((x_a)^{1/2}-(x_a)^{-1/2})\cdots
&((q_\al^{-k_\al} x_a)^{1/2}-(q_\al^{-k_\al} x_a)^{-1/2}).
}
$$
Then $\widetilde{\mu}^*=\widetilde{\mu}$ and
$$
\langle p_b\widetilde{\mu}\tga^{-1}\rangle\ =\
(\langle p_b^*\,\widetilde{\mu}\tga\rangle)^*.
$$
  
Let 
$$ 
\eqalignno{
&\widetilde{\mu}=\widetilde{\mu}^{(k)},\ 
\widetilde{\mu}'=\widetilde{\mu}^{(k+v)},\ b'=b-r_v, \cr
& p=p^{(k)}_b,\ p' = p^{(k+v)}_{b-r_v} = 
 (g^{k}(b))^{-1}\g_v^{(k)}(p).
}
$$ 
The notations are from Proposition \ref\GP.

\proclaim { Key Lemma}
$$
\eqalignno{
&\langle (p')^*\, \widetilde{\mu}' \tga \rangle\ =\ 
q^{(r_k,r_k)/2-(b-r_{k},b-r_{k})/2} (d_{k+v}/d_k)  
&\eqnu\cr
\label\key\eqnum*
&\times \prod_{\al\in R_+,\nu_\al\in v} 
\Bigl
(q_\al^{(b-r_k,\al)/2+k_\al}-
q_\al^{(r_k-b,\al)/2}\Bigr)
\langle (p)^*\, \widetilde{\mu} \tga \rangle,\cr
 &d_k= 
|W(r_k)|^{-1}\prod_{\al\in R_+,k_\al\neq 0} 
\Bigl( { q_\al^{((r_k,\al)+k_\al)/2}-
q_\al^{-((r_k,\al)+k_\al)/2}
\over
 q_\al^{(r_k,\al)/2}-
q_\al^{-(r_k,\al)/2} }
 \Bigr).
}
$$
\label\KEY\theoremnum*
\endproclaim 
{\it Proof.}
One has:
$$
\eqalign{
&\langle  (p')^*\,\widetilde{\mu}'\tga\rangle =
(d_{k+v}/d_k)\langle (\x_v^{(k)})^2 (p')^*\,\widetilde{\mu}\tga\rangle\cr
&=(d_{k+v}/d_k)\langle (-1)^{\ka_v}(\x_v^{(k)}\tga) (\x_v^{(k)}p')^*
\widetilde{\mu}\rangle,
}
\eqnu
\label\redmu\eqnum*
$$
where $\ka_v$ is the number of roots $\al>0$ with $\nu_\al\in v$.
See Lemma 4.4 and formula (5.7) from [C2].
We will treat the factors $ \x_v^{(k)}$ in two different ways.

First,
$$
\eqalign{
&\x_v^{(k)}\tga=q^{(r_k,r_k)/2-(b-r_k,b-r_k)/2} \y_v^{(k)}(\tga)
}
\eqnu
$$
thanks to formula (\ref\epy) applied to $\x_v^{(k)}$, that
is a linear combination of $e_c^{(k)}$ for $c\in W(r_v)$.
Indeed, $\y_v^{(k)}(m_{c})$ is proportional to $\x_v^{(k)}$
(or zero) for all $B_-\ni c\succeq -r_v$. Applying any $Y_a$ to $e_b$
we get  linear combinations of the $e$-polynomials from the
same orbit. 

Second,
$$ 
\eqalignno{
&\x_v^{(k)}p'= (g^{k}(b))^{-1}\y_v^{(k)}(p).
}
$$ 
Substituting, let us combine two $\y$ together using 
the unitarity of  $\{Y\}$
with respect to the pairing $\langle f g^*\,\widetilde{\mu}\rangle$:  
$$
\eqalign{
&\langle (p')^*\,\widetilde{\mu}'\tga\rangle=
(d_{k+v}/d_k)\langle \tga\ (\{\y_v^{(k)}\}^2 (p))^*\,
\widetilde{\mu}\rangle.
}
\eqnu
\label\redmu\eqnum*
$$

Because the constant term is $W$-invariant, we may multiply 
 $\y^2$ on the left  by (const)$\p_{r_v}^{(k)}$ for the $t$-symmetrization
from (\ref\symmetr) with the idempotent normalization. 
As to $\y^2$, it can be replaced by
$(-1)^{\ka_v}\y\bar{\y}$
for
$$
\bar{\y}_v^{(k)}  = \prod_{\al\in R_+,\nu_\al\in v}
((t_\al Y_{a})^{1/2}-
(t_\al Y_{a})^{-1/2}),\ a=\al^\vee
$$
(see [C2], formula (5.10)).
Now we use formula (\ref\Lf). Collecting all the factors
together we get the required.
\proofbox

\comment
Let us take any set $k=\{k_{\nu_1}\ge k_{\nu_2}\}\in \Z_+$ and put
$$
\eqalign{
t(k)=\{ q^{2k_\nu/\nu}\},\   k\cdot r=\sum_\nu k_\nu r_\nu,\ 
p_b^{(k)}= p_b^{t(k)}.
}
\eqnu
$$ 
\endcomment

The remaining  part of the calculation is based on
 the following chain of the shift operators, that will be applied
to $p^{(0)}_{-r_k}=m_{-r_k}$ 
one after another:
$$
\eqalignno{
&\g_{\nu_R}^{(k-e)}\g_{\nu_R}^{(k-2e)}\cdots 
\g_{\nu_R}^{(k-se)}\g_{\nu_1}^{(k-se-v)}
\cdots \g_{\nu_1}^{(0)},
&\eqnu
}
$$
where 
$ k_{\nu_1}=s+d,\ k_{\nu_2}=s,\ 
v=\{v_\nu\},\  v_{\nu_1}=1,\ v_{\nu_2}=0,$
the set $\{1,1\}$ is denoted by $e$. We assume that $d\ge 0$
choosing the components properly. So
$\ k-se= dv$.

The product
$$
q^{(r_{k-e},r_{k-e})/2-(r_{k},r_{k})/2} (d_{k}/d_{k-e})\cdots
q^{-(r_{v},r_{v})/2} (d_{v}/d_0)
$$
is equal to $q^{-(r_{k},r_{k})/2}d_{k}$.
Taking into consideration that $p_{-r_k}(1)=|W(r_k)|$,
we come to the final formula for 
$\langle \widetilde{\mu}^{(k)}\tga\rangle$.
Conjugating,
$$
\eqalignno{
&\langle \tga^{-1}\widetilde{\mu}^{(k)}\rangle\  =\ 
&\eqnu\cr
\label\tmuga\eqnum*
&q^{(r_k,r_k)}\prod_{\al\in R_+}
\prod_{1\le j\le k_\al}
\Bigl( q_\al^{(r_k,\al)/2+1/2}-q_\al^{-(r_k,\al)/2+1/2-j}
\Bigr).
}
$$
Multiplying by $\mu/\widetilde{\mu}$, we verify
(\ref\mehtamu) for $k_\al\in \Z_+$.

If $k$ are arbitrary we can still use the key lemma.
It gives that the ratio $\phi(t)$ of the left and
the right  hand sides of (\ref\mehtamu) is a periodic function
of $t_\nu$ with respect to the shift $t_\nu\to t_\nu q_\nu$.
We may replace the constant term by a contour integral 
in $x$ for $|t_\al q_\al|<|x_a|<|t_\al^{-1}|,\ $ 
$a=\al^\vee, \al>0$, provided that
$|t_\nu|, |q_\nu|<1$. If $t$ has only one component, then
$\phi(t)$ is analytic and $q$-periodic for such $q,t$
and has to be $1$.
Otherwise, $t=(t_{\nu_1},t_{\nu_2})$ and we conclude that $\phi(t)$ 
depends only on $\tau=t_{\nu_1}^{\nu_1/2}t_{\nu_2}^{-\nu_1/2}$. Taking
$t_{\nu_2}=1$ and applying the shift operator for
$v=\nu_1$, we see that $\phi(\tau)=\phi(q\tau)$ and $\phi=1$
(cf. [O1]). The calculation is completed.

%
%
%
\section{ Jackson pairing}

The formulas from the previous sections can be generalized
for Jackson integrals taken instead of the constant term.
We fix $\xi\in \C^n$ and set $\langle f\rangle_\xi\equal
|W|^{-1}\sum_{w\in W, b\in B}f(q^{w(\xi)+b})$ following
[C1]. Here $f$ is a Laurent polynomial or any function well-defined
at $\{q^{w(\xi)+b}\}$ provided the convergence. We 
assume that $|q|<1$. One may also operate in formal series
in terms of non-negative powers of $q$ considering
$q^{(\xi,b_i)}$ as independent parameters.

Recall that  $x_a(q^\xi)=q^{(a,\xi)}, \ga(q^z)=q^{(z,z)/2}$.
For instance, 
$$\langle \ga\rangle_\xi= \sum_{a\in B} q^{(\xi+a,\xi+a)/2}=
\tga^{-1}(q^\xi)q^{(\xi,\xi)/2}.$$ 

Applying the shift operators for $k_\nu\in \Z_+$ we readily establish that
$$
\eqalignno{
&\langle \ga{\mu}^\circ\rangle_\xi\  =\ 
 |W|^{-1}q^{(r_k,r_k)/2}\langle \ga p_{r_k}\rangle_\xi \cr
&\times \prod_{\al\in R_+}
\prod_{1\le j\le k_\al}
\Bigl( 1-q_\al^{-(r_k,\al)/2-j}
\Bigr), \for \cr
&\mu^\circ \equal \prod_{\al\in R_+}\prod_{0\le j\le k_\al-1}
\Bigl(1-q_\al^{-j}X_a^{-1})(1-q_\al^{-j-1}X_a)
\Bigr).
&\eqnu\cr
\label\mucirck\eqnum*
}
$$
We use the same Lemma \ref\KEY.

The analytic continuation to any $t$ is 
$$
\eqalign{
&\mu^\circ\  =\  \prod_{\al \in R_+}
\prod_{i=0}^\infty {(1-x_a t_\al^{-1}q_\al^{i}) 
(1-x_a^{-1}t_\al^{-1}q_\al^{i+1})
\over
(1-x_a q_a^{i}) (1-x_a^{-1}q_\al^{i+1})}.
}
\eqnu
\label\mucirc\eqnum*
$$
In fact it is $\mu^{-1}(x,t^{-1})$.
We assume in the next theorem that $\mu^\circ(q^{w(\xi)+b})$ is
well-defined, i.e. $(\al,\xi)\not\in \Z$ for all $\al\in R_+$.

\proclaim{ Theorem}
Given $b,c\in B$ and the corresponding polynomials $\ep_b,\ep_c$,
$$
\eqalignno{
&\langle
\ep_b\, \ep_c^*\,\ga{\mu^\circ}
\rangle_\xi\ =\
 q^{-(b_\#,b_\#)/2-(c_\#,c_\#)/2 +(r_k,r_k)} 
\ep_c(q^{b_\#})\langle \ga{\mu^\circ}\rangle_\xi,
&\eqnu \cr
\label\eejack\eqnum*
& \langle \ga{\mu^\circ}\rangle_\xi\ =\ 
|W|^{-1}\langle \ga\rangle_\xi 
 \prod_{\al\in R_+}\prod_{ j=0}^{\infty}\Bigl({ 
1- t_\al^{-1}q_\al^{-(r_k,\al)+j}\over
1- q_\al^{-(r_k,\al)+j} }\Bigr). 
&\eqnu
\label\memujack\eqnum*  
}
$$
\label\EEJACK\theoremnum*
\endproclaim

In this formulas $t$ is arbitrary provided the existence
of $\ep_{b,c}$.  The right hand side of (\ref\memujack) is 
replaced by the limit when $k_\nu\in \Z_+$. We use the values
at these points and the shift operators to deduce (\ref\memujack)
from (\ref\mucirck) (see the end of the previous section).
There is one more similar formula corresponding to the 
conjugation of (\ref\epep). These formulas can be extended
even to the special cases when $\mu^\circ (q^{w(\xi)+b})$ has
singularities, for instance, to  spherical representations 
of double affine Hecke algebras from [C1]. One has  to renormalize
$\mu$ as follows. 

Let us 
switch from $\mu^\circ$
to 
$\widehat{\mu}\equal \mu^\circ/\mu^\circ(q^{-r_k})$.
See [C1], Proposition 4.2 and
[M5].
Setting $\La(bw)=R_+^a\cap (bw)^{-1}(R_-^a)$,  
$$
\eqalignno{
& \widehat{\mu}(q^{w(\xi)+b}) =\widehat{\mu}(q^{w(\xi)+b})^*=
\prod_{[\al,j]\in \La(bw)}
\Bigl(
{
t_\al^{-1/2}-t_\al^{1/2} q_\al^{(\al,\xi)+j}\over
t_\al^{1/2}-t_\al^{-1/2} q_\al^{(\al,\xi)+j}
}
\Bigr),
&\eqnu
\label\muval\eqnum*
}
$$
where the conjugation of
$\xi$ is trivial:
$(q^{(\xi,a)})^*=q^{-(\xi,a)}$.
We note that (\ref\muval) is the same for $\mu$ 
instead of $\mu^\circ$.

Let $\xi=-r_k$ for generic $k$. Then 
$\widehat{\mu}(q^{w(\xi)+b})$ is always well-defined
and non-zero only for $\pi_b=b\om_b^{-1}$ [C1]. 
The previous theorem in this  case
reads as follows:

\proclaim{ Theorem}
For $\xi=-r_k,\  b, c\in B$,
$$
\eqalignno{
&\langle
\ep_b\, \ep_c^*\,\ga\widehat{\mu}
\rangle_\xi =
 q^{-(b_\#,b_\#)/2-(c_\#,c_\#)/2 +(r_k,r_k)} 
\ep_c(q^{b_\#})\langle \ga\widehat{\mu}\rangle_\xi,
&\eqnu \cr
\label\hatjack\eqnum*
& \langle \ga\widehat{\mu}\rangle_\xi\ =\ 
|W|^{-1}\langle \ga\rangle_\xi 
\prod_{\al\in R_+}\prod_{ j=1}^{\infty}\Bigl({ 
1- q_\al^{(r_k,\al)+j}\over
1-t_\al^{-1}q_\al^{(r_k,\al)+j} }\Bigr). 
&\eqnu
\label\hatmu\eqnum*  
}
$$
\label\HATMU\theoremnum*
\endproclaim

\proclaim {Corollary}
For the function $\e_q$ from
(\ref\exelas), arbitrary $b,c\in B$,
and the constant $\langle\ga\rangle_{r_k}=
\ga(q^{r_k})\tga^{-1}(q^{r_k}),$ 
$$
\eqalignno{
&\e_q(q^{b_\#}, q^{c_\#})\langle\ga\rangle_{r_k}^2 \ =\ 
\ep_c(q^{b_\#})|W| \langle \ga\widehat{\mu}\rangle_{-r_k},
&\eqnu \cr
\label\ebc\eqnum*
& \e_q(x, q^{c_\#})\langle\ga\rangle_{r_k}\ =\ 
\ep_c(x) 
\prod_{\al\in R_+}\prod_{ j=1}^{\infty}\Bigl({ 
1- q_\al^{(r_k,\al)+j}\over
1-t_\al^{-1}q_\al^{(r_k,\al)+j} }\Bigr). 
&\eqnu
\label\hatmux\eqnum*
}
$$
\label\HATMUX\theoremnum*
\endproclaim
{\it Proof.}
One has:
$$
\eqalignno{
&\e_{q}(q^{b_\#},q^{c_\#}))\tga^{-1}(q^{b_\#})\tga^{-1}(q^{c_\#})\ =\cr 
&\sum_{a\in B} q^{(a_\#,a_\#)/2 -(r_k,r_k)} 
\ {\ep_a(q^{b_\#})\ \ep_a^*(q^{c_\#})\over
\langle \ep_a,\ep_a\rangle_1 }=\cr
&\sum_{a\in B} q^{(a_\#,a_\#)/2 -(r_k,r_k)} 
\ \ep_b(q^{a_\#})\ \ep_c^*(q^{a_\#})
\widehat{\mu}(q^{a_\#})=\cr
&q^{-(r_k,r_k)}|W|\langle \ep_b\, \ep_c^*\,\ga\widehat{\mu}
\rangle_{-r_k}=\cr
&q^{-(r_k,r_k)} q^{-(b_\#,b_\#)/2-(c_\#,c_\#)/2 +(r_k,r_k)} 
\ep_c(q^{b_\#})|W|\langle \ga\widehat{\mu}\rangle_{-r_k}.
&\eqnu
\label\ebcjack\eqnum*
}
$$
Here we used the duality and the formula
$\widehat{\mu}(q^{b_\#})=\langle \ep_b,\ep_b\rangle_1^{-1}$
from Theorem 5.6 [C1].
The function $q^{(b_\#,b_\#)/2}\tga^{-1}(q^{b_\#})$ does not
depend on $b\in B$, so we get the constant $\langle\ga\rangle_{r_k}^2$.

The second formula (\ref\hatmux) is an analytic continuation of the first.
The ratio $\psi=\e_q(x, q^{c_\#})/\ep_c(x)$ is $B$-periodic and
with poles in the $B$-periodic set of zeros of $\tga^{-1}(x)$.
Therefore $\psi$ is analitic for all $x_i\neq 0$ and must be a
constant. 
\proofbox

Thus $\e_q$ can be introduced more conceptually as a unique
extension of the non-symmetric polynomials in the class of
meromorphic functions of non-zero $\{x_i,\la_j\}$ with the poles
in the zeros of the Gaussians $\tga^{-1}(x),\tga^{-1}(\la)$.
It must satisfy (\ref\yelet) and be invariant with respect to the
action of the intertwining operators. 

Let us give a symmetric version of the above results (see Section 1). 
The reproducing kernel is given by the formula in terms of
symmetric Macdonald's polynomials:
$$
\eqalignno{
& q^{(r_k,r_k)/2}\p(x,\la)\tga_x^{-1}\tga_\la^{-1}\ =\
\sum_{b\in B_-} q^{(b,b)/2 -(r_k,b)} 
\ {p_b(x)\ p_b(\la^{-1})\langle \De\rangle\over
\langle p_b(x) p_b(x^{-1})\De\rangle}.
&\eqnu\cr
\label\pexla\eqnum*
}
$$
The right hand side is analytic for  non-zero $x_i,\la_j$
when $|q|<1$. The function $\p$ is 
$x\leftrightarrow \la$ symmetric and satisfies
the equations
$$
\eqalignno{
&f(Y_1,\ldots,Y_n)(\p(x,\la))\ =\ f(\la^{-1})\p(x,\la) \for f\in \C[x]^W.
&\eqnu\cr
\label\symeqpi\eqnum*
}
$$
 
It extends $\{p_c, c\in B_-\}$ with the same coefficient
of proportionality as in  (\ref\hatmux):
$$
\eqalignno{
& \p(x, q^{c-r_k})\langle\ga\rangle_{r_k} = 
p_c(x)(p_c(q^{r_k}))^{-1} 
\prod_{\al\in R_+}\prod_{ j=1}^{\infty}\Bigl({ 
1- q_\al^{(r_k,\al)+j}\over
1-t_\al^{-1}q_\al^{(r_k,\al)+j} }\Bigr). 
&\eqnu
\label\hatmux\eqnum*
}
$$
In fact it is a corollary of  Theorem \ref\PLANJACK for $\xi=-r_k$.

We note that the calculation of the 
constant $\langle \ga{\mu^\circ}\rangle_\xi$
from Theorem \ref\EEJACK and that for $\widehat{\mu}$ instead of $\mu^\circ$
is directly related to the Aomoto conjecture [A,Ito] proved
in [M5]. Actually it gives another way to check this
conjecture.

%
%
%
%
%
%
\AuthorRefNames [BGG]
\references
\vfil

[A]
\name {K. Aomoto},
{On product formula for Jackson integrals associated with root
systems}, Preprint (1994).

[AW]
\name{R. Askey}, and \name{J. Wilson},
{  Some basic hypergeometric orthogonal polynomials
that generalize Jacobi polynomials},
Memoirs AMS {  319} (1985).

[B]
\name{N. Bourbaki},
{ Groupes et alg\`ebres de Lie}, Ch. {\bf 4--6},
Hermann, Paris (1969).

[C1]
\name{I. Cherednik},
{ Intertwining operators of double affine Hecke algebras},
Selecta Math. (1997).

[C2]
\bibline,{ Double affine Hecke algebras and  Macdonald's
conjectures},
Annals of Mathematics {141} (1995), 191-216.

[C3]
\bibline, 
{ Macdonald's evaluation conjectures and
difference Fourier transform},
Inventiones Math. {122} (1995),119--145.

[C4]
\bibline, 
{ Nonsymmetric Macdonald polynomials },
IMRN {10} (1995), 483--515.

[DS]
\name{J.F. van Diejen}, and \name{J.V. Stockman},
{ Multivariable $q$-Racah polynomials},
Preprint CRM-2369 (1996). 

[D]
\name{C.F. Dunkl}, 
{ Hankel transforms associated to finite reflection groups},
Contemp. Math. {138} (1992), 123--138.

[J]
\name{M.F.E. de Jeu},
{The Dunkl transform }, Invent. Math. {113} (1993), 147--162.

[H]
\name {S. Helgason},
{Groups and geometric analysis}, 
Academic Press, New York (1984).

[Ito]
\name{M. Ito}, {On a theta product formula
for Jackson integrals associated with root
systems of rank two},
Preprint (1996).

[KL]
\name{D. Kazhdan}, and \name{ G. Lusztig},
{  Proof of the Deligne-Langlands conjecture for Hecke algebras},
Invent.Math. {  87}(1987), 153--215.

[KK]
\name{B. Kostant}, and \name{ S. Kumar},
{  T-Equivariant K-theory of generalized flag varieties,}
J. Diff. Geometry{  32}(1990), 549--603.

[M1]
\name{I.G. Macdonald}, {  Some conjectures for root systems},
SIAM J.Math. Anal. { 13}:6 (1982), 988--1007.

[M2]
\bibline, {  A new class of symmetric functions },
Publ.I.R.M.A., Strasbourg, Actes 20-e Seminaire Lotharingen,
(1988), 131--171 .

[M3]
\bibline, {  Orthogonal polynomials associated with root 
systems},Preprint(1988).

[M4]
\bibline, { Affine Hecke algebras and orthogonal polynomials},
S\'eminaire Bourbaki { 47}:797 (1995), 01--18.

[M5]
\bibline, { A formal identity for affine root systems},
Preprint (1996).

[O1]
\name{E.M. Opdam}, 
{  Some applications of hypergeometric shift
operators}, Invent.Math.{  98} (1989), 1--18.

[O2]
\bibline, { Harmonic analysis for certain representations of
graded Hecke algebras}, 
Preprint Math. Inst. Univ. Leiden W93-18 (1993).

[O3]
\bibline, { Dunkl operators, Bessel functions and the discriminant of
a finite Coxeter group},
Composito Mathematica {85} (1993), 333--373.

\endreferences

\bye